\documentclass[10pt]{amsart}

\usepackage{amsmath}
\usepackage{amssymb}
\usepackage{bm}
\usepackage{graphicx}
\usepackage{xcolor}
\usepackage{hyperref}
\hypersetup{hidelinks}
\usepackage{url}
\usepackage{algpseudocode}
\usepackage{fancyhdr}
\usepackage{mathtools}
\usepackage{tikz-cd}
\usepackage{xy}
\xyoption{all}
\usepackage{stmaryrd}
\usepackage{calrsfs}

\theoremstyle{plain}
\newtheorem{theorem}{Theorem}[section]
\newtheorem{lemma}[theorem]{Lemma}

\newtheorem{prop}[theorem]{Proposition}
\newtheorem{corollary}[theorem]{Corollary}

\newtheorem{conjecture}[theorem]{Conjecture}

\theoremstyle{definition}

\newtheorem{example}[theorem]{Example}

\newtheorem{prob}[theorem]{Problem}

\theoremstyle{remark}
\newtheorem{remark}[theorem]{Remark}

\numberwithin{equation}{section}
\theoremstyle{definition}
\newtheorem*{example*}{Proposed Example}

\newcommand{\cc}{\mathbb{C}}

\newcommand{\nn}{\mathbb{N}}
\newcommand{\pp}{\mathbb{P}}

\newcommand{\qq}{\mathbb{Q}}
\newcommand{\rr}{\mathbb{R}}
\newcommand{\zz}{\mathbb{Z}}

\DeclareMathOperator{\mcd}{mcd}
\DeclareMathOperator{\CD}{CD}
\providecommand\ldb{\llbracket}
\providecommand\rdb{\rrbracket}

\newcommand{\gp}{\mathcal{G}}

\newcommand{\num}{\mathsf{n}}
\newcommand{\den}{\mathsf{d}}

\definecolor{reviewred}{RGB}{180,0,0}
\definecolor{reviewgreen}{RGB}{0,110,45}
\definecolor{reviewblue}{RGB}{0,0,0}
\definecolor{reviewgray}{RGB}{95,95,95}

\newcommand{\backgroundadd}[1]{\textcolor{reviewblue}{#1}}

\title{Factorization in rational monogenic semidomains II}

\author{Felix Gotti}
\address{Department of Mathematics\\MIT\\Cambridge, MA 02139}
\email{fgotti@mit.edu}

\author{Keshav Karumbu}
\address{Jesuit High School\\Portland, OR 97225}
\email{keshavkarumbu0@gmail.com}

\author{Hyun-Jin Kim}
\address{Liberal Arts and Science Academy\\Austin, TX 78721}
\email{joshuahyunjink@gmail.com}

\author{Vincent Wang}
\address{Liberal Arts and Science Academy\\Austin, TX 78721}
\email{vincentwangtao@gmail.com}

\date{October 5, 2026}



\begin{document}

\begin{abstract}
	For $q \in \qq_{ > 0}$, let $S_q := \nn_0[q]$. We investigate divisibility and factorization in the multiplicative monoid $S_q^*$. We characterize the AP, PSP, D, GL, GCD, and quasi-GCD properties, and record the equivalence of the FF, BF, and ACCP properties using earlier results. We prove that $S_q^*$ is Furstenberg if and only if its technical monoid $T_q$ is Furstenberg. We construct explicit families of non-Furstenberg and non-atomic rational monogenic semidomains and classify strong atoms in the atomic cases. For $q < 1/2$ outside the ACCP cases, we characterize atoms and ACCP elements through unique polynomial representations and prove that the MCD and quasi-ACCP properties are equivalent. We also obtain MCD-finiteness criteria and construct pairs with arbitrarily many nonassociate MCDs in every nonfactorial case. When $\num(q) > 1$, the denominator $\den(q)$ has at least two distinct prime divisors, and every prime-power factor of $\den(q)$ exceeds $\num(q)$, we determine the atoms of $T_q$ explicitly: they are the primes dividing $\den(q)$ that are smaller than $\num(q)$. Consequently, $S_q^*$ is not Furstenberg, and it is MCD-finite but neither MCD nor quasi-ACCP. Finally, in the non-ACCP range, we characterize ACCP elements by a uniform bound on the coefficients of their polynomial representations.
\end{abstract}

\maketitle


\bigskip
\section{Introduction}
\label{sec:intro}

The arithmetic of a semiring can differ substantially from that of the ring it generates. For a positive rational number $q$, the semiring
\[
S_q := \nn_0[q] = \{f(q) : f \in \nn_0[x]\}
\]
provides a natural setting in which to investigate this difference. Its elements admit polynomial representations with nonnegative integer coefficients, but subtraction is unavailable. Consequently, divisibility in its multiplicative monoid $S_q^* := S_q \setminus \{0\}$ depends on the existence of such representations. Understanding how polynomial representations govern multiplicative factorizations is the main concern of this paper.
\smallskip

The present paper is a continuation of the work of Deng, Gotti, and Zeng~\cite{DGZ26} on rational monogenic semidomains. They characterized the UF, HF, FF, and Krull properties in this class, proved that the BF property is equivalent to the ascending chain condition on principal ideals (ACCP), and introduced the divisor-closed monoid $T_q$ of technical fractions. Its elements are the elements of $S_q^*$ whose reduced numerators have all their prime divisors among those of $\den(q)$. Their investigation left open whether every rational monogenic semidomain is atomic. We give explicit negative answers to this question and continue their study through the Furstenberg property, maximal common divisors, and the ACCP for individual elements.
\smallskip

The additive arithmetic of monogenic semidomains provides an earlier foundation for this investigation. Gotti and Gotti~\cite{GG18} briefly determined the positive rationals whose generated semidomains have atomic additive monoids. Gotti~\cite{fG19} proved that increasing positive monoids of ordered fields have the FF property, a result that applies to the additive monoid of $\nn_0[\alpha]$ for every real number $\alpha > 1$. Chapman, Gotti, and Gotti~\cite{CGG20} studied sets of lengths, elasticity, and other factorization invariants of the additive monoids $(S_q, +)$. Correa-Morris and Gotti~\cite{CG22} extended the study to positive algebraic generators, characterizing additive atomicity and factoriality and proving the equivalence of the additive ACCP, BF, and FF properties.
\smallskip

More recent work develops this additive theory in several directions. Chen, Gotti, Lu, and Yao~\cite{CGLY26} studied algebraic monogenic semidomains from a valuation-theoretic perspective, characterizing the absence of additive atoms and identifying parameters for which the additive monoid is a valuation monoid. Chen, Gotti, and Yao~\cite{CGY26} extended the setting to finitely generated semidomains, investigating the effect of adjoining generators and proving the equivalence of additive bounded and finite factorization for positive real generators. Related results on additive factorizations in evaluation semirings and Laurent evaluation semirings were obtained by Ajran et al.~\cite{ABLST23} and Zhu~\cite{sZ22}. Dani et al.~\cite{DDGLPVZ26} investigated additive atoms and strong atoms in cyclic semidomains.
\smallskip

The simultaneous study of additive and multiplicative factorization was developed by Baeth, Chapman, and Gotti~\cite{BCG21}, who investigated bi-atomic classes of positive semirings. Dani, Gotti, Li, and Paladiya~\cite{DGLP26} studied algebraic monogenic semidomains, including their unit groups, criteria for non-atomicity, and finite factorization properties. Gotti, Graia, Han, and Liang~\cite{GGHL26} established the Bi-UF Positive Conjecture for quadratic monogenic semirings and studied related half-factoriality questions. Gotti, Graia, and Panse~\cite{GGP26} subsequently settled the Bi-LF Positive Conjecture for semidomains and classified rational and transcendental monogenic positive semirings satisfying prescribed quasi-factoriality conditions. These results place rational monogenic semidomains within a broader study of the interaction between the two monoid structures of a semidomain.
\smallskip

Polynomial factorization is also relevant because the evaluation map transfers polynomial identities to $S_q$. Campanini and Facchini~\cite{CF19} studied factorization in $\nn_0[x]$, while Gotti and Polo~\cite{GP23a,GP23b} investigated arithmetic and weaker forms of atomicity in polynomial semidomains. Gonzalez, Gotti, and Panpaliya~\cite{GGP25} studied the ascent of almost atomicity and quasi-atomicity to polynomial extensions and monoid semidomains. Baeth and Gotti~\cite{BG20} considered cyclic Puiseux information semialgebras in their study of factorization in upper triangular matrices, and Cesarz, Chapman, McAdam, and Schaeffer~\cite{CCMS09} studied elasticity in semirings arising from positive systems.
\smallskip

The divisibility properties considered here have also been studied in more general monoids. The GCD, quasi-GCD, MCD, and MCD-finite properties appear in work of Blitz, Gotti, Han, and Liang~\cite{BGHL26} on monoid algebras. Liang, Wang, and Zhong~\cite{LWZ24} investigated the existence of MCDs in Puiseux monoids and its connections with atomicity and the ACCP. Related work on the Furstenberg property appears in papers by Lebowitz-Lockard~\cite{nL19}, Lin, Rabinovitz, and Zhang~\cite{LRZ23}, and Liu, Rodr\'iguez, and Tirador~\cite{LRT24}. Dani et al.~\cite{DGHLS25} study the quasi-ACCP, while Angerm\"uller~\cite{gA22}, Fadinger et al.~\cite{FFNSW26}, and Gotti and Krause~\cite{GK24} study strong atoms. We recall the definitions used here in Section~\ref{sec:background}. For background on Puiseux monoids and factorization theory, see Chapman, Gotti, and Gotti~\cite{CGG21}, Geroldinger and Halter-Koch~\cite{GH06}, and Anderson and Gotti~\cite{AG22}.
\smallskip

After the algebraic preparations in Section~\ref{sec:algebra}, we prove in Section~\ref{sec:uf-gl-d-ap} that the UF, PSP, D, GL, and AP properties coincide in this class and hold precisely when $q \in \nn \cup \nn^{-1}$. In Section~\ref{sec:bf-ff}, we record that the FF and BF properties and the ACCP are equivalent, holding precisely when $q \ge 1$, $q \in \nn^{-1}$, or $\den(q)$ is a prime power. This classification follows by combining the technical-monoid results and factorization theorems of Deng, Gotti, and Zeng~\cite{DGZ26}. We also strengthen their approximation result for technical fractions.
\smallskip

In Section~\ref{sec:furstenberg}, we prove that $S_q^*$ is Furstenberg if and only if $T_q$ is Furstenberg. Outside the ACCP cases, finiteness of the atom set of $T_q$ forces failure of the Furstenberg property. In particular, suppose that $\num(q) > 1$, that $\den(q)$ has at least two distinct prime divisors, and that $p^{\nu_p(\den(q))} > \num(q)$ for every prime $p$ dividing $\den(q)$. Then the atoms of $T_q$ are exactly the primes dividing $\den(q)$ that are smaller than $\num(q)$, so $S_q^*$ is not Furstenberg and hence is not atomic. We obtain further families of non-atomic rational monogenic semidomains and determine the strong atoms whenever $S_q^*$ is atomic.
\smallskip

In Section~\ref{sec:gcd-mcd}, we show that the GCD and quasi-GCD properties have the same parameter classification as the UF property. We establish MCD-finiteness under the preceding prime-power hypothesis and construct pairs with arbitrarily many nonassociate MCDs in every nonfactorial case. In Section~\ref{sec:accp}, we prove that, when $q < 1/2$ and $S_q^*$ does not satisfy the ACCP, an element satisfies the ACCP if and only if it has a unique polynomial representation; in this range, the MCD and quasi-ACCP properties are equivalent. Finally, throughout the non-ACCP range, we characterize ACCP elements by a uniform bound on the coefficients of all their polynomial representations. This criterion also applies to elements with infinitely many representations.


\bigskip
\section{Background}
\label{sec:background}

In this section, we fix the notation and terminology used throughout the paper.


\medskip
\subsection{General Notation}

We let $\zz$, $\qq$, $\rr$, and $\cc$ denote the sets of integers, rational numbers, real numbers, and complex numbers, respectively. We let $\pp$, $\nn$, and $\nn_0$ denote the sets of primes, positive integers, and nonnegative integers, respectively. We set $\nn^{-1} := \{1/n : n \in \nn\}$ and call its elements \emph{unit fractions}. For $r$ and $s$ in $\rr$, set $\ldb r, s \rdb := \{n \in \zz : r \le n \le s\}$. For $S \subseteq \rr$ and $r \in \rr$, set $S_{ \ge r} := \{s \in S : s \ge r\}$ and define $S_{ > r}$ analogously. We use the convention $\inf\varnothing := + \infty$. An integer $b \in \nn$ is a \emph{prime power} if $b = p^e$ for some $p \in \pp$ and $e \in \nn$. For $b \in \nn$ with $b > 1$ and a prime $p$ dividing $b$, the \emph{full prime-power factor} of $b$ corresponding to $p$ is $p^{\nu_p(b)}$, where the valuation is defined below. In this paper, ``every prime-power factor of $b$'' refers to these full factors. For $r \in \qq \setminus \{0\}$, write $r = \num(r)/\den(r)$ in lowest terms, with $\den(r) > 0$. For $p \in \pp$, first set
\[
\nu_p(n) := \max\{k \in \nn_0 : p^k \mid n\}
\]
for each $n \in \zz \setminus \{0\}$. The \emph{$p$-adic valuation} is the map $\nu_p \colon \qq \setminus \{0\} \to \zz$ defined by
\[
\nu_p(r) := \nu_p(\num(r)) - \nu_p(\den(r)).
\]
For the zero rational, we use the conventions $\num(0) := 0$ and $\den(0) := 1$. For positive integers $d$ and $b$, the notation $d \mid b^\infty$ means that $d \mid b^N$ for some $N \in \nn_0$; equivalently, every prime divisor of $d$ divides $b$. We say that a positive rational is \emph{supported on} a set of primes if every prime divisor of its numerator or denominator lies in that set. Congruences of rational numbers modulo an integer are interpreted by inverting denominators relatively prime to that integer. For $n \in \nn$, the \emph{Euler totient function} $\phi(n)$ counts the integers in $\ldb 1, n \rdb$ relatively prime to $n$.
\smallskip

For a nonzero polynomial \backgroundadd{$f \in R[x]$ over a commutative ring or semiring $R$} and \backgroundadd{$j \in \nn_0$}, the symbols $\deg f$ and $[x^j]f$ denote the degree and the coefficient of $x^j$, respectively. Coefficients beyond the degree are taken to be zero. \backgroundadd{We also set $[x^j]0 := 0$ for every $j \in \nn_0$.} A polynomial in $\zz[x]$ is \emph{primitive} if the greatest common divisor of its coefficients is $1$. A nonzero polynomial is \emph{monic} if its leading coefficient is $1$, and is \emph{squarefree} over $\qq$ if no nonconstant polynomial square divides it in $\qq[x]$. For \backgroundadd{$n \in \nn$}, the \emph{$n$th cyclotomic polynomial} $\Phi_n$ is the monic polynomial whose roots are the complex roots of unity of exact order $n$; it is irreducible in $\qq[x]$. We use Gauss's lemma in its usual polynomial form: a primitive polynomial in $\zz[x]$ dividing an integer polynomial in $\qq[x]$ also divides it in $\zz[x]$. This is distinct from the monoid GL property defined below.
\smallskip

For a commutative ring or semiring $R$, the notation $R[x]$ denotes the polynomials in one indeterminate $x$ with coefficients in $R$. A coefficient sequence $(c_i)_{i \ge 0}$ in $R$ is \emph{eventually zero} if there exists $N \in \nn_0$ such that $c_i = 0$ for every $i \ge N$; thus, writing a polynomial as $\sum_{i \ge 0} c_i x^i$ involves only finitely many nonzero terms. A \emph{monomial} is a polynomial $cx^j$ with $c \in R \setminus \{0\}$ and $j \in \nn_0$. The \emph{leading coefficient} of a nonzero polynomial $f \in R[x]$ is $[x^{\deg f}]f$, and its \emph{constant coefficient} is $[x^0]f$. If $f \in R[x]$ and $\rho$ lies in a commutative ring or semiring containing $R$ then \emph{evaluation at $\rho$} is the map sending $f = \sum_{i \ge 0} c_i x^i$ to $f(\rho) := \sum_{i \ge 0} c_i\rho^i$.
\smallskip

Let $X_0 \in \rr$, and let $F$ and $G$ be positive real-valued functions defined on $(X_0, + \infty)$. Then the notation $F(X) = O(G(X))$ means that $F(X) \le C G(X)$ for some constant $C \in \rr_{ > 0}$ and all sufficiently large $X \in \rr$; the notation $F(X) = o(G(X))$ means that $F(X)/G(X) \to 0$. We write $F(X) \sim G(X)$ when their ratio tends to $1$.



\medskip
\subsection{Commutative Monoids}

Throughout this paper, a \emph{monoid} is a commutative cancellative semigroup with an identity, usually written multiplicatively. Throughout this paper, we tacitly assume that every monoid homomorphism preserves identity elements and every submonoid inherits the identity element. 
\smallskip

Let $M$ be a monoid. A \emph{unit} is an invertible element, and $M^\times$ denotes the group consisting of all the units of~$M$. The monoid $M$ is called \emph{reduced} if $M^\times = \{1\}$.  The \emph{reduced monoid} of $M$ is the monoid consisting of all the cosets of $M$ by $M^\times$:
\[
	M_{\mathrm{red}} := \{aM^\times : a \in M\}
\]
under the natural multiplication induced by $M$. Additive submonoids of $\qq_{ \ge 0}$ are reduced monoids, which are referred to as \emph{Puiseux monoids}. A Puiseux monoid that is cofinite in $\nn_0$ is called a \emph{numerical monoid}. The \emph{Frobenius number} of a numerical monoid $M_0$, denoted by $F(M_0)$, is the maximum of $\nn_0 \setminus M_0$ if $M_0 \neq \nn_0$ and $-1$ otherwise. For relatively prime integers $a$ and $b$ with $1 < a < b$ and $n \in \nn$, the geometric Frobenius formula of Ong and Ponomarenko~\cite{OP08} is
\begin{equation} \label{eq:geometric-frobenius}
	F\bigl(\langle a^n, a^{n - 1}b, \dots, b^n\rangle\bigr) = \frac{(a - 1)b^{n + 1} - (b - 1)a^{n + 1}}{b - a}.
\end{equation}
The same expression is symmetric in $a$ and $b$ and therefore applies when $a > b > 1$ as well.
\smallskip

The \emph{Grothendieck group} of $M$, denoted by $\gp(M)$, is an abelian group equipped with a monoid embedding $\iota \colon M \to \gp(M)$ such that $\iota(M)$ generates $\gp(M)$ as a group. For a subset $S$ of $M$, we let $\langle S \rangle$ denote the smallest submonoid of $M$ containing $S$, and we call $\langle S \rangle$ the submonoid of $M$ \emph{generated by}~$S$. A monoid is \emph{finitely generated} if it can be generated by a finite set.

\medskip
\subsection{Atomicity and Factorization}

An element $a \in M \setminus M^\times$ is called an \emph{atom} (or an \emph{irreducible element}) if, for all $b$ and $c$ in $M$, the equality $a = bc$ ensures that either $b$ or $c$ is a unit. Write $\mathcal A(M)$ for the set of atoms. An element of a monoid is called \emph{atomic} if it is a unit or can be written as a finite product of atoms. We say that a monoid is \emph{atomic} if every element is atomic. We say that a monoid is \emph{antimatter} if it contains no atoms at all. An atom $a$ is a \emph{strong atom} if $a^n$ has a unique factorization for every $n \in \nn$. We write $\mathcal A_{\mathrm{str}}(M)$ for the set of strong atoms. Strong atoms are also called absolutely irreducible elements; see Geroldinger and Halter-Koch~\cite{GH06}, Angerm\"uller~\cite{gA22}, Fadinger et al.~\cite{FFNSW26}, and Gotti and Krause~\cite{GK24}.
\smallskip

Let $\mathsf{Z}(M)$ denote the free commutative monoid on the set $\mathcal{A}(M_{\text{red}})$, and let $\pi \colon \mathsf{Z}(M) \to M_{\text{red}}$ be the unique monoid homomorphism fixing the set $\mathcal{A}(M_{\text{red}})$. For each $b \in M_{\text{red}}$, set $\mathsf{Z}(b) := \pi^{-1}(b)$ and call the elements of $\mathsf{Z}(b)$ \emph{factorizations} of~$b$ in~$M$. An element $b \in M$ is called \emph{factorial} if $|\mathsf{Z}(b)| = 1$. The monoid $M$ is a \emph{unique factorization monoid} (UFM) if $|\mathsf{Z}(b)| = 1$ for all $b \in M$, which means that it consists of factorial elements. The monoid $M$ is said to have the \emph{finite factorization} (FF) property if every element of~$M$ has a nonempty finite set of factorizations; that is, $1 \le |\mathsf{Z}(b)| < \infty$ for all $b \in M$. By definition, every UFM is an FFM. Given a factorization $z \in \mathsf{Z}(M)$, the \emph{length} of $z$, denoted by $|z|$, is the number of atoms of $z$ (counting repetitions). For each $b \in M_{\text{red}}$, set
\[
	\mathsf{L}(b) := \{ |z| : z \in \mathsf{Z}(b) \}.
\]
An atomic monoid is said to have the \emph{bounded factorization} (BF) property if each element has a bounded set of factorization lengths. By definition, every monoid having the FF property also has the BF property. For background on the BF and FF properties, see Anderson and Gotti~\cite{AG22}.

\medskip
\subsection{Divisibility}

For $b,c \in M$, we write $b \mid_M c$ and say that $b$ \emph{divides} $c$ in $M$ provided that $c = ab$ for some $a \in M$. A submonoid $M_0 \subseteq M$ is \emph{divisor-closed} if every divisor in $M$ of an element of $M_0$ belongs to $M_0$, in which case, $M_0^\times = M_0 \cap M^\times$. Divisibility induces a partial order $\preceq$ on $M_{\mathrm{red}}$ as follows: for $a$ and $b$ in $M$, we have $aM^\times \preceq bM^\times$ precisely when $a \mid_M b$. If $M$ is reduced then this is a partial order on $M$ itself. 
\smallskip

A nonunit $p \in M$ is called \emph{prime} if whenever $p \mid_M bc$ for some $b,c \in M$ we can conclude that either $p \mid_M b$ or $p \mid_M c$. A monoid is factorial precisely when every nonunit is a product of primes. Notice that every prime is an atom, and a monoid is said to have the \emph{AP property} provided that all its atoms are prime. Every UFM is an AP monoid. A monoid is called \emph{Furstenberg} if every nonunit is divisible by an atom. Observe that every factorial monoid is atomic while every atomic monoid is Furstenberg.
\smallskip

For a nonempty subset $E \subseteq M$, a \emph{common divisor} is an element dividing every member of $E$. Write $\CD_M(E)$ for the set of common divisors. A common divisor divisible by every common divisor of $E$ is a \emph{greatest common divisor} (GCD). A common divisor $d \in M$ is a \emph{maximal common divisor} (MCD) if the only common divisors of $E/d := \{e/d : e \in E\}$ are units. Write $\gcd_M(E)$ and $\mcd_M(E)$ for the respective sets. Equivalently, MCDs are the maximal elements of $\CD_M(E)$ under divisibility modulo associates; a \emph{maximal element} has no strictly larger element in this order.
\smallskip

The monoid is a \emph{GCD monoid} or an \emph{MCD monoid} if every finite nonempty subset has a GCD or an MCD, respectively. It is \emph{quasi-GCD} (q-GCD) if every such subset has at most one MCD up to associates, and \emph{MCD-finite} if every such subset has finitely many MCDs up to associates, allowing none. Every GCD is an MCD. Thus, GCD implies both MCD and q-GCD, while q-GCD implies MCD-finite. The q-GCD and MCD-finite properties in this monoid-theoretic form are studied by Blitz, Gotti, Han, and Liang~\cite{BGHL26}, and the MCD property of Puiseux monoids is studied by Liang, Wang, and Zhong~\cite{LWZ24}.
\smallskip

Elements $u$ and $v$ of $M$ are \emph{relatively prime} if all their common divisors are units. \backgroundadd{For $u$, $v$, and $w$ in $M$,} the \emph{D property} requires that $u \mid_M vw$ and relative primality of $u$ and $v$ imply $u \mid_M w$. The \emph{Gauss Lemma (GL) property} requires that relative primality of $u$ to each of $v$ and $w$ imply relative primality of $u$ to $vw$. The PSP, GL, and AP properties originate in the theory of integral domains, where it is classical that PSP implies GL and GL implies AP; see Anderson and Zafrullah~\cite{AZ07} and Gotti and Zafrullah~\cite{GZ23}. A closely related monoid version of the GL property appears in the work of Blitz, Gotti, Han, and Liang~\cite{BGHL26}. The D property is the monoid counterpart of the divisibility condition through which PSP implies GL; we are not aware of a standard name for it in the literature.

\medskip
\subsection{Ideals and Chain Conditions}

Let \backgroundadd{$M$ be a monoid}. A nonempty subset $I \subseteq M$ is an \emph{ideal} if $MI \subseteq I$. It is \emph{principal} if $I = aM$ for some $a \in M$, and \emph{finitely generated} if $I = a_1M \cup \cdots \cup a_kM$ for \backgroundadd{some $k \in \nn$ and} $a_1, \dots, a_k \in M$. An element \backgroundadd{$a \in M$} satisfies the \emph{ascending chain condition on principal ideals} (ACCP) if every ascending chain $aM = a_0M \subseteq a_1M \subseteq \cdots$ eventually stabilizes. \backgroundadd{An ascending chain $(I_n)_{n \ge 0}$ of ideals \emph{stabilizes} if there exists $N \in \nn_0$ such that $I_n = I_N$ for every $n \ge N$. An \emph{ACCP element} is an element satisfying the ACCP.} The monoid satisfies the ACCP if every element does. A divisor of an ACCP element satisfies the ACCP, and every ACCP element is atomic: otherwise repeated splitting of a non-atomic nonunit gives an infinite strictly ascending chain of principal ideals.
\smallskip

The monoid satisfies the \emph{quasi-ACCP} if, for every finite nonempty subset $E \subseteq M$, there are a common divisor $d \in M$ of $E$ and an element $e \in E$ such that $e/d$ satisfies the ACCP. The monoid satisfies the \emph{almost ACCP} if it is atomic and satisfies the quasi-ACCP. Equivalently, the almost ACCP is obtained by requiring the common divisor $d$ in the quasi-ACCP condition to be atomic: applying this condition to a singleton proves atomicity of the monoid. These conditions were introduced by Gotti and Li~\cite{GL22}; we follow the terminology of Dani et al.~\cite{DGHLS25}. The implications
\[
\text{UF} \ \Longrightarrow \ \text{FF} \ \Longrightarrow \ \text{BF} \ \Longrightarrow \ \text{ACCP} \ \Longrightarrow \ \text{quasi-ACCP} \ \Longrightarrow \ \text{MCD}
\]
hold for commutative cancellative monoids. To justify the last implication, choose $d$ and $e$ as above. If $E/d$ has no MCD then successively enlarging common divisors $c_0 = 1, c_1, c_2, \dots$ in $M$ yields a strictly ascending chain of principal ideals generated by $(e/d)/c_i$ for $i \in \nn_0$, starting at $(e/d)M$ and contradicting its ACCP. Multiplying an MCD of $E/d$ by $d$ gives an MCD of $E$.
\smallskip

For subsets $E$ and $F$ of $\gp(M)$, set $(F:E) := \{g \in \gp(M) : gE \subseteq F\}$. For a nonempty ideal $I \subseteq M$, put $I^{-1} := (M:I)$ and $I_v := (M:(M:I))$. The ideal is \emph{divisorial}, or a \emph{$v$-ideal}, if $I = I_v$. The monoid is \emph{$v$-Noetherian} if every ascending chain of divisorial ideals stabilizes. Principal ideals are divisorial, so $v$-Noetherianity implies the ACCP. A finitely generated ideal is \emph{primitive} if it is contained in no proper principal ideal and \emph{super-primitive} if $I^{-1} = M$. The monoid has the \emph{primitive super-primitive} (PSP) property if every primitive ideal is super-primitive. Every super-primitive ideal is primitive: containment in $aM$ for a nonunit $a \in M$ would give $a^{-1} \in I^{-1} \setminus M$.

\medskip
\subsection{Semirings and Rational Monogenic Semidomains}

A \emph{commutative semiring} is a set $S$ with associative and commutative addition and multiplication, identities $0$ and $1$, distributivity, and $0s = 0$ for every $s \in S$. The multiplication need not be cancellative on a general semiring. A \emph{subsemiring} contains both identities and is closed under both operations. A \emph{semiring homomorphism} preserves both operations and both identities, and an \emph{isomorphism} is a bijective homomorphism. A \emph{semidomain} $S$ is a subsemiring of a field. Its additive monoid is $(S, +)$, and its multiplicative monoid is $S^* := S \setminus \{0\}$. We write $S^\times := (S^*)^\times$. All factorization and divisibility properties ascribed to $S$ refer to $S^*$. In particular, $a \mid_S b$ for nonzero elements $a$ and $b$ of $S$ means $a \mid_{S^*} b$. \backgroundadd{A \emph{unique factorization semidomain} (UFS) is a semidomain $S$ for which $S^*$ is factorial. We also use \emph{FF-monoid} for a monoid with the FF property.}
\smallskip

An \emph{integral domain} is a commutative ring with $1 \ne 0$ and no zero divisors; it is a semidomain after embedding in its field of fractions. \backgroundadd{For an integral domain $R$, its \emph{field of fractions} consists of the fractions $a/b$ with $a$ and $b$ in $R$ and $b \ne 0$. For $a$, $b$, $c$, and $d$ in $R$ with $b \ne 0$ and $d \ne 0$, the equality $a/b = c/d$ holds precisely when $ad = bc$. The field $\qq(x)$ is the field of fractions of $\qq[x]$ and is called the \emph{rational function field} in the indeterminate $x$ over $\qq$.} A \emph{semifield} is a semidomain whose nonzero elements are units; examples are $\qq_{ \ge 0}$ and $\rr_{ \ge 0}$. The polynomial semiring $\nn_0[x]$ is a semidomain inside $\qq(x)$. Its multiplicative monoid satisfies the FF property. Indeed, a divisor of a nonzero polynomial $f \in \nn_0[x]$ has degree at most $\deg f$ and every coefficient at most the largest coefficient of $f$, so there are finitely many divisors. Moreover, in a factorization $f = gh$ with $g$ and $h$ in $\nn_0[x]^*$ and $h \ne 1$, the positive integer $\deg g + g(1)$ is strictly smaller than $\deg f + f(1)$, which proves atomicity.
\smallskip

For $\rho \in \cc$, the \emph{monogenic semidomain generated by $\rho$} is
\[
S_\rho := \nn_0[\rho] = \{f(\rho) : f \in \nn_0[x]\}.
\]
A semidomain is \emph{rational monogenic} if it is isomorphic to $S_q$ for some $q \in \qq$. We restrict throughout the results to $q > 0$. If $q = 0$ then $S_q = \nn_0$. For \backgroundadd{$a$ and $b$ in $\nn$}, if $q = -a/b < 0$ is in lowest terms then $-1 = (a - 1) + bq \in S_q$, so $S_q = \zz[q] = \zz[1/b]$; hence that case reduces to a localization of the integers. Here $\zz[1/b] := \{z/b^k : z \in \zz, \ k \in \nn_0\}$.
\smallskip

Fix $q = a/b \in \qq_{ > 0}$ in lowest terms, \backgroundadd{where $a$ and $b$ belong to $\nn$}. Every element of $S_q$ has denominator dividing $b^\infty$. A \emph{polynomial representation} of $r \in S_q$ is a polynomial $f \in \nn_0[x]$ satisfying $f(q) = r$. Define the homomorphism $\omega_q \colon S_q^* \longrightarrow (\nn_0, +)$ by
\[
\omega_q(r) := \sum_{\substack{p \in \pp \\ p \nmid b}} \nu_p(r),
\]
and call $T_q := \omega_q^{-1}(0)$ the \emph{technical monoid}. The sum is finite and nonnegative because the reduced denominator of $r$ is supported on the prime divisors of $b$. For $0 < q < 1$ with $a > 1$, define the real number
\[
B_q := \frac{a - 1}{1 - q}.
\]
For integers $a$ and $b$ with $1 < a < b$, define the \emph{digit-numerator set}
\[
\mathcal D_{a, b} := \left\{\sum_{i = 0}^n c_i a^i b^{n - i} :
n \in \nn_0, \ c_0, \dots, c_n \in \ldb 0, a - 1 \rdb, \ c_n > 0\right\}.
\]
\backgroundadd{An element of $T_q$ is a \emph{technical fraction}, and a \emph{technical atom} is an atom of $S_q^*$ that belongs to $T_q$.} The monoid $T_q$, threshold $B_q$, and set $\mathcal D_{a, b}$ will control technical divisibility, coefficient growth, and MCD-finiteness, respectively.
\smallskip

\backgroundadd{For the related settings mentioned in the introduction, a complex number $\alpha$ is \emph{algebraic} if it is a root of a nonzero polynomial in $\qq[x]$. A \emph{Laurent polynomial} over a commutative ring or semiring $R$ is a finite sum of terms $c_i x^i$ with $c_i \in R$ and $i \in \zz$. If $M$ is an additive commutative monoid and $R$ is a commutative ring then the \emph{monoid algebra} $R[M]$ consists of formal finite sums $\sum_{m \in M} c_m X^m$ with $c_m \in R$, with coefficientwise addition and multiplication determined by $X^mX^n = X^{m + n}$ for $m$ and $n$ in $M$. A \emph{lattice monoid} is an additive submonoid of $\zz^d$ for some $d \in \nn$.}

\bigskip
\section{Algebraic Considerations}
\label{sec:algebra}

We recall the basic properties of the technical monoid and establish a criterion for unique polynomial representations. The resulting divisor correspondence will be used throughout the paper.

\medskip
\subsection{The Technical Monoid}

Fix $q \in \qq_{ > 0}$. Recall the technical monoid $T_q = \omega_q^{-1}(0)$ from Section~\ref{sec:background}. Additivity and nonnegativity of $\omega_q$ imply that $T_q$ is divisor-closed. Equivalently,
\[
T_q = \{r \in S_q^* : \num(r) \mid \den(q)^\infty\}.
\]
The membership condition $r \in S_q^*$ is essential: not every positive rational supported on denominator primes belongs to $S_q^*$.

If $\num(q) > 1$ then every polynomial representation of an element of $T_q$ has positive constant coefficient. Indeed, a representation with zero constant coefficient gives divisibility by $q$, whereas $\omega_q(q) > 0$. Therefore $T_q \subseteq \qq_{ \ge 1}$. Every unit of $S_q^*$ belongs to $T_q$, as does its inverse. Both are at least $1$, and their product is $1$, so $S_q^*$ is reduced. If $q = 1/b$ for some $b \in \nn$ then $S_q = \{N/b^k : N, k \in \nn_0\}$, and its units are exactly the positive rationals supported on prime divisors of $b$. In this case $T_q = (S_q^*)^\times$ and $S_q^*$ is factorial, with atoms represented by primes not dividing $b$. For positive integral $q$, the equality $S_q = \nn_0$ gives factoriality as well.

\medskip
\subsection{Polynomial Representations}

We first characterize unique polynomial representations for every positive rational parameter.

\begin{lemma} \label{lem:unique-polynomial-representation}
	Let $a$ and $b$ in $\nn$ be relatively prime, set $q := a/b \in \qq_{>0}$, and let $c \in S_q$. The following conditions are equivalent.
	\begin{enumerate}
		\item[(a)] $c$ has exactly one polynomial representation.
		\smallskip
		
		\item[(b)] $c$ has a representation $f = \sum_{i \ge 0} c_i x^i \in \nn_0[x]$, where $(c_i)_{i \ge 0}$ is an eventually zero sequence in $\nn_0$, such that $c_0 < a$ and $c_i < \min\{a, b\}$ for every $i \ge 1$.
	\end{enumerate}
	In particular, when $0 < q < 1$, uniqueness is equivalent to the existence of a representation whose coefficients are all less than $a$.
\end{lemma}

\begin{proof}
	(a) $\Rightarrow$ (b): Choose a representation $f = \sum_{i \ge 0} c_i x^i \in \nn_0[x]$ of $c$, where $(c_i)_{i \ge 0}$ is an eventually zero sequence in $\nn_0$. Fix $i \in \nn_0$. If a coefficient $c_i$ is at least $a$ then replacing $a x^i$ by $b x^{i + 1}$ gives a different representation. If $i \ge 1$ and $c_i \ge b$ then replacing $b x^i$ by $a x^{i - 1}$ gives a different representation.
	\smallskip
	
	(b) $\Rightarrow$ (a): Conversely, suppose that $f$ satisfies~(b), and let $g \in \nn_0[x]$ be another representation of $c$. Gauss's lemma gives $g - f = (bx - a)h$ for some $h \in \zz[x]$. Assume $h \ne 0$, set $d := \deg h \in \nn_0$, write $h = \sum_{i = 0}^d h_i x^i$ with $h_0, \dots, h_d \in \zz$, and set $h_{-1} := 0$. If $h_d < 0$ then the coefficient of $x^{d + 1}$ in $g$ is at most $(b - 1) - b < 0$. Hence $h_d > 0$. Let $j \in \ldb 0, d \rdb$ be the smallest index with $h_j > 0$. Then $h_{j - 1} \le 0$, and the coefficient of $x^j$ in $g$ is
	\[
	c_j + b h_{j - 1} - a h_j \le (a - 1) - a < 0,
	\]
	which is a contradiction. Therefore $h = 0$ and $g = f$.
\end{proof}

A unique polynomial representation transfers the divisibility structure from the polynomial semiring to its value.

\begin{prop} \label{prop:unique-divisors}
	Let $q \in \qq_{ > 0}$, and suppose $r \in S_q^*$ has the unique polynomial representation $f \in \nn_0[x]$. Evaluation at $q$ induces an isomorphism between the divisibility-ordered set of divisors of $f$ in $\nn_0[x]^*$ and the divisibility-ordered set of divisors of $r$ in $S_q^*$. Every divisor of $r$ has a unique representation, and $r$ has finitely many divisors and finitely many factorizations into atoms. In particular, $r$ satisfies the ACCP.
\end{prop}

\begin{proof}
	Let $u$ and $v$ in $S_q^*$ satisfy $r = uv$, and choose representations $g$ and $h$ in $\nn_0[x]$ of $u$ and $v$, respectively. Uniqueness gives $f = gh$. If $g' \in \nn_0[x]$ also represents $u$ then $g'h = f = gh$, so cancellation in $\zz[x]$ gives $g' = g$. Thus, all divisors have unique representations. The same argument proves injectivity of evaluation on polynomial divisors; surjectivity follows from the chosen $g$ and $h$. Moreover, if $g(q)$ divides $h(q)$ and $h$ divides $f$, multiplying a representation of $h(q)/g(q)$ by $g$ and by $f/h$ shows that $g$ divides $h$. Hence the correspondence preserves divisibility in both directions. The polynomial divisor set is finite, as proved in the background, and its factorizations into atoms correspond to those of $r$. Every ascending chain of principal ideals starting at $rS_q^*$ uses divisors of $r$, so it stabilizes.
\end{proof}

\bigskip
\section{UF, PSP, GL, D, and AP Properties}
\label{sec:uf-gl-d-ap}

Fix $q \in \qq_{ > 0}$. We prove that every atom of $S_q^*$ is prime if and only if $q \in \nn \cup \nn^{-1}$, and use this characterization to classify the UF, PSP, D, and GL properties.

\begin{lemma} \label{lem:atom-condition-sq}
	Let $q = m/n \in \qq_{ > 0}$ with coprime integers $m \ge 2$ and $n \ge 2$, and let $N \in \nn_0[x]$ be irreducible. If every coefficient of $N$ is less than $\min\{m, n\}$ then $N(q)$ is an atom of $S_q^*$.
\end{lemma}

\begin{proof}
	Lemma~\ref{lem:unique-polynomial-representation} shows that $N$ is the unique representation of $N(q)$. Proposition~\ref{prop:unique-divisors} transfers irreducibility of $N$ to $N(q)$.
\end{proof}

We now use Lemma~\ref{lem:atom-condition-sq} to show that $S_q^*$ is not AP whenever $q \in \qq_{ > 0} \setminus(\nn \cup \nn^{-1})$.

\begin{theorem} \label{thm:non-ap}
	If $q \in \qq_{ > 0} \setminus(\nn \cup \nn^{-1})$ then $S_q^*$ is not AP.
\end{theorem}

\begin{proof}
	Write $q = m/n$ in lowest terms, where $m \in \nn_{\ge 2}$ and $n \in \nn_{\ge 2}$, so $q \ne 1$. Define the elements $a$, $b$, and $c$ of $S_q^*$ by
	\[
	a := 1 + q^3, \qquad b := 1 + q, \qquad c := 1 + q^2 + q^4.
	\]
	The identity
	\[
	bc = (1 + q)(1 + q^2 + q^4) = (1 + q^3)(1 + q + q^2)
	\]
	shows that $a \mid_{S_q^*}bc$.
	
	Each of $1 + x^3$, $1 + x$, and $1 + x^2 + x^4$ is irreducible in $\nn_0[x]$. Indeed, its value at $1$ is a prime integer, so in every factorization one factor has value $1$ at $1$ and must therefore be a monomial with coefficient~$1$. Since the constant term of the product is $1$, that factor must be the constant polynomial~$1$. Every coefficient of these polynomials is at most $1 < \min\{m, n\}$, so Lemma~\ref{lem:atom-condition-sq} shows that $a$, $b$, and $c$ are atoms of $S_q^*$. Moreover, $a \ne b$ because $q \ne 1$, and
	\[
	c - a = q^2(q^2 - q + 1) > 0.
	\]
	Since $1$ is the only unit of $S_q^*$, $a$ is associate to neither $b$ nor $c$. An atom can divide another atom only if they are associates. Thus, $a \nmid_{S_q^*}b$ and $a \nmid_{S_q^*}c$, although $a \mid_{S_q^*}bc$. Therefore $a$ is not prime, and $S_q^*$ is not AP.
\end{proof}

For completeness, we prove the standard implications among these divisibility properties directly for commutative cancellative monoids.

\begin{prop}
	Let $M$ be a commutative cancellative monoid. The following implications hold for $M$:
	\begin{equation} \label{eq:divisibility-implications}
		\begin{tikzcd}
			{\textbf{ UF }}
			\arrow[r, Rightarrow]
			&
			{\textbf{ PSP }}
			\arrow[r, Rightarrow]
			&
			{\textbf{ D }}
			\arrow[r, Rightarrow]
			&
			{\textbf{ GL }}
			\arrow[r, Rightarrow]
			&
			{\textbf{ AP }}.
		\end{tikzcd}
	\end{equation}
\end{prop}

\begin{proof}
	Suppose first that $M$ is a UFM. For a primitive finitely generated ideal $I \subseteq M$ and an element $z \in \gp(M)$, prime exponents show that $zI \subseteq M$ forces $z \in M$: a negative prime exponent in $z$ would force that prime to divide every element of $I$. Hence UF implies PSP.
	
	To argue that PSP implies the D property, let $u \in M$, $v \in M$, and $w \in M$, suppose that $u$ and $v$ are relatively prime, and assume that $u \mid_M vw$. The ideal $I := uM \cup vM$ is primitive. Since $(w/u)I \subseteq M$, we have $w/u \in I^{-1} = M$, proving $u \mid_M w$. Here $I^{-1} := (M:I)$.
	
	To verify that the D property implies the GL property, let $u \in M$, $v \in M$, and $w \in M$, suppose that $u$ is relatively prime to both $v$ and $w$, and let $d \in M$ divide both $u$ and $vw$. Then $d$ is relatively prime to $v$, so D gives $d \mid_M w$. Relative primality of $u$ and $w$ forces $d$ to be a unit.
	
	Finally, suppose that $M$ is GL, let $v \in M$ and $w \in M$, and suppose that an atom $u \in \mathcal A(M)$ divides $vw$ but divides neither $v$ nor $w$. Every nonunit divisor of $u$ is associate to $u$, so $u$ is relatively prime to both $v$ and $w$. GL would make $u$ relatively prime to $vw$, contradicting $u \mid_M vw$. Thus, GL implies AP.
\end{proof}

Factoriality for $q \in \nn \cup \nn^{-1}$, Theorem~\ref{thm:non-ap}, and the implication chain~\eqref{eq:divisibility-implications} yield the following classification.

\begin{corollary} \label{cor:uf-ap}
	Let $q \in \qq_{ > 0}$. The following conditions are equivalent.
	\begin{enumerate}
		\item[(a)] $S_q$ is a UFS;
		\smallskip
		
		\item[(b)] $S_q$ has the PSP property;
		\smallskip
		
		\item[(c)] $S_q$ has the D property;
		\smallskip
		
		\item[(d)] $S_q$ has the GL property;
		\smallskip
		
		\item[(e)] $S_q$ has the AP property;
		\smallskip
		
		\item[(f)] $q \in \nn \cup \nn^{-1}$.
	\end{enumerate}
\end{corollary}

\bigskip
\section{Bounded and Finite Factorization Properties}
\label{sec:bf-ff}

We record the equivalence of the FF and BF properties and the ACCP for rational monogenic semidomains. This follows by combining Proposition~4.9(2) and Theorems~5.7--5.8 of Deng, Gotti, and Zeng~\cite{DGZ26}. The approximation lemma below strengthens the statement that $1$ is a limit point of $T_q$ by requiring $(z - 1)/q^s \in S_q$ for a prescribed $s \in \nn_0$. This additional divisibility will be used to construct common technical divisors.


\begin{lemma} \label{lem:tech-approximation}
	Let $a \in \nn$ and $b \in \nn$ be relatively prime integers with $2 \le a < b$, and set $q := a/b \in \qq \cap (0, 1)$. Suppose that $b$ has at least two distinct prime divisors. For every $s \in \nn_0$ and $\epsilon \in \rr_{ > 0}$, there exists a nonunit $z \in T_q$ such that
	\[
	\frac{z - 1}{q^s} \in S_q,
	\]
	with $1 < z < 1 + \epsilon$.
\end{lemma}

\begin{proof}
	Write $b = mn$, where $m \ge 2$ and $n \ge 2$ are coprime integers; such a decomposition exists because $b$ has at least two distinct prime divisors. Choose $L \in \nn_0$ with $L \ge s$ sufficiently large that
	\[
	(a - 1)bq^L < \epsilon.
	\]
	The numbers $\log m$ and $\log n$ are linearly independent over $\qq$: a rational relation would imply equality between positive powers of the coprime integers $m > 1$ and $n > 1$. Hence the set of differences $U\log m - V\log n$ with $U \in \nn$ and $V \in \nn$ is dense in $\rr$, even when both indices are required to exceed a prescribed bound (see Deng, Gotti, and Zeng~\cite[Lemma~4.8]{DGZ26}). We can therefore choose $U \in \nn$ and $V \in \nn$, with $\phi(a^L)V > L$, such that
	\[
	\log\bigl(1 + (a - 1)bq^L\bigr) < \phi(a^L)(U\log m - V\log n) < \log(1 + \epsilon).
	\]
	Set $u := \phi(a^L)U \in \nn$ and $v := \phi(a^L)V \in \nn$, so that $v - L \ge 1$. Exponentiating yields
	\[
	1 + (a - 1)bq^L < \frac{m^u}{n^v} < 1 + \epsilon.
	\]
	Since $\gcd(a, m) = \gcd(a, n) = 1$, Euler's theorem gives
	\[
	m^u \equiv n^v \equiv 1 \pmod{a^L},
	\]
	so $m^u - n^v = Ta^L$ for some $T \in \nn$. The lower bound on $m^u/n^v$ gives
	\[
	\frac{Ta^L}{n^v} > (a - 1)b\frac{a^L}{b^L}.
	\]
	Multiplication by $m^v n^v/a^L$ yields
	\[
	Tm^v > (a - 1)b^{v - L + 1}.
	\]
	Let $F \in \zz$ be the Frobenius number of the additive numerical monoid
	\[
	\langle a^{v - L}, a^{v - L - 1}b, \dots, b^{v - L}\rangle.
	\]
	By the theorem of Ong and Ponomarenko~\cite{OP08}, this number satisfies
	\begin{align*}
		F & = b^{v - L - 1}(ab - a - b)
		+ \frac{(b - 1)a^2(b^{v - L - 1} - a^{v - L - 1})}{b - a}\\
		& = \frac{(a - 1)b^{v - L + 1} - (b - 1)a^{v - L + 1}}{b - a}\\
		& < (a - 1)b^{v - L + 1} < Tm^v.
	\end{align*}
	Thus, there are $d_0, \dots, d_{v - L} \in \nn_0$ such that
	\[
	Tm^v = \sum_{i = 0}^{v - L}d_i a^i b^{v - L - i}.
	\]
	Set $z := m^u/n^v \in \qq_{ > 0}$. Division by $b^{v - L}$ gives the stronger conclusion
	\[
	\frac{z - 1}{q^L} = \frac{Tm^v}{b^{v - L}} = \sum_{i = 0}^{v - L}d_iq^i \in S_q.
	\]
	In particular, $z = 1 + q^L\sum_{i = 0}^{v - L}d_iq^i \in S_q$. Both the numerator and the denominator of $z$ have all their prime divisors among those of $b$, so $z \in T_q$. Since $L \ge s$, we obtain
	\[
	\frac{z - 1}{q^s} = q^{L - s}\cdot \frac{z - 1}{q^L} \in S_q.
	\]
	Moreover, since $a \ge 2$, Lemma~\ref{lem:unique-polynomial-representation} shows that $1$ has the unique polynomial representation $1$, and hence $S_q^*$ is reduced. Therefore $z > 1$ is a nonunit, and by construction, $1 < z < 1 + \epsilon$.
\end{proof}

The approximation lemma and the earlier factorization criteria give the following classification.


\begin{theorem} \label{thm:ff-bf-accp}
	For $q \in \qq_{ > 0}$, let $S_q$ be the monogenic semidomain generated by $q$. Then the following conditions are equivalent.
	\begin{enumerate}
		\item[(a)] $S_q^*$ satisfies the FF property.
		\smallskip
		
		\item[(b)] $S_q^*$ satisfies the BF property.
		\smallskip
		
		\item[(c)] $S_q^*$ satisfies the ACCP.
		\smallskip
		
		\item[(d)] Either $q \in \nn^{-1}$ or $1$ is not a limit point of $T_q$.
		\smallskip
		
		\item[(e)] Either $q \in [1, + \infty) \cup \nn^{-1}$ or $\den(q)$ is a prime power.
	\end{enumerate}
\end{theorem}

\begin{proof}
	(a) $\Rightarrow$ (b): This follows directly from the corresponding definitions.
	\smallskip
	
	(b) $\Rightarrow$ (e): Assume that $S_q^*$ satisfies the BF property and suppose, by way of contradiction, that~(e) fails. Then $q = a/b$ in lowest terms, where $a \in \nn$ and $b \in \nn$ satisfy $2 \le a < b$ and $b$ has at least two distinct prime divisors. Since $q$ is not a unit fraction, Section~\ref{sec:algebra} gives $T_q \subseteq \qq_{ \ge 1}$, and Deng, Gotti, and Zeng~\cite[Theorem~5.8]{DGZ26} give $\inf(T_q \setminus \{1\}) > 1$.
	
	On the other hand, Lemma~\ref{lem:tech-approximation} applied with $s = 0$ (or, alternatively, the result of Deng, Gotti, and Zeng~\cite[Proposition~4.9(2)]{DGZ26}) gives for every $\epsilon \in \rr_{ > 0}$ a nonunit $z \in T_q$ such that $1 < z < 1 + \epsilon$. Hence $\inf(T_q \setminus \{1\}) \le 1$, a contradiction.
	\smallskip
	
	(e) $\Rightarrow$ (a): If $q = 1$ then $S_q = \nn_0$, which is factorial. If $q > 1$, set $\delta := \min\{2, q\} \in \rr_{ > 1}$. Every nonunit of $S_q^*$ is at least $\delta$. For each fixed $r \in S_q^*$, the set $S_q \cap [1, r]$ is finite. Indeed, a polynomial representation $\sum_{i \ge 0} e_iq^i \le r$, with an eventually zero sequence $(e_i)_{i \ge 0}$ in $\nn_0$, has $e_i = 0$ whenever $q^i > r$, and $0 \le e_i \le r/q^i$ for every remaining index $i \in \nn_0$. Every product of nonunits equal to $r$ has length at most $\log r/\log\delta$, so repeated nontrivial factorization terminates in atoms. Every factor belongs to the finite set $S_q \cap [1, r]$, giving only finitely many factorizations of $r$. Thus, $S_q^*$ satisfies the FF property. If $q^{-1} \in \nn$ then $S_q^*$ is factorial and hence satisfies the FF property. The remaining prime-power denominator case follows from Deng, Gotti, and Zeng~\cite[Theorem~5.7]{DGZ26}.
	\smallskip
	
	(b) $\Leftrightarrow$ (c) $\Leftrightarrow$ (d): If $q \in \nn^{-1}$ then $S_q^*$ is factorial, so it satisfies both the BF property and the ACCP, while condition~(d) holds by definition. If $q \ge 1$ then~(a) holds by the argument above, so~(b) and~(c) hold; moreover,~(d) holds because every nonunit of $S_q^*$ is at least $\delta > 1$ when $q > 1$, while $T_q = \{1\}$ when $q = 1$. Hence we may assume that $q \in (0, 1) \setminus \nn^{-1}$. Then $T_q \subseteq \qq_{ \ge 1}$, and Deng, Gotti, and Zeng~\cite[Theorem~5.8]{DGZ26} give
	\[
	S_q^*\text{ satisfies BF}
	\quad\Longleftrightarrow\quad
	S_q^*\text{ satisfies ACCP}
	\quad\Longleftrightarrow\quad
	\inf(T_q \setminus \{1\}) > 1.
	\]
	Here $\inf\varnothing := + \infty$. Since $T_q \subseteq \qq_{ \ge 1}$, the final inequality is equivalent to $1$ not being a limit point of $T_q$, as required.
\end{proof}

\bigskip
\section{The Furstenberg Property}
\label{sec:furstenberg}

For each $q \in \qq_{ > 0}$, we study the Furstenberg property of $S_q^*$. This property is related to the infinitude of primes and irreducibles; see Clark~\cite{pC17}, Lebowitz-Lockard~\cite{nL19}, and Gotti and Krause~\cite{GK26}. For the Furstenberg property in Puiseux monoids and lattice monoids, see Lin, Rabinovitz, and Zhang~\cite{LRZ23} and Liu, Rodr\'iguez, and Tirador~\cite{LRT24}, respectively.

Recall that $T_q$ is divisor-closed in $S_q^*$. In particular,
\[
T_q^\times = T_q \cap (S_q^*)^\times,
\qquad
\mathcal A(T_q) = T_q \cap \mathcal A(S_q^*).
\]
These identities allow us to pass between atoms of $T_q$ and technical atoms of $S_q^*$.

\medskip
\subsection{Set of Atoms}

Throughout this subsection, unless otherwise stated, let $a \in \nn$ and $b \in \nn$ satisfy $2 \le a < b$ and $\gcd(a, b) = 1$, and set $q := a/b \in \qq$. Assume that $b$ has at least two distinct prime divisors. By Theorem~\ref{thm:ff-bf-accp}, these are precisely the parameters for which $S_q^*$ fails the ACCP. In this case $(S_q^*)^\times = \{1\}$ and $T_q \subseteq \qq_{ \ge 1}$: every element of $T_q$ has a polynomial representation with positive constant coefficient, since otherwise it is divisible by $q$, contrary to $\omega_q(q) > 0$ and $\omega_q(T_q) = \{0\}$. Thus, a unit and its inverse, both in $T_q$, must equal~$1$.

Recall the threshold used to study $\mathcal A(S_q^*)$:
\[
B_q := \frac{\num(q) - 1}{1 - q} = \frac{\den(q)(\num(q) - 1)}{\den(q) - \num(q)}.
\]

This is the constant denoted by $B_{q, \infty}$ by Deng, Gotti, and Zeng~\cite[Corollary~3.4]{DGZ26}. The next lemma proves the scaled membership statement needed here.


\begin{lemma} \label{lem:technical-division}
	Let $a \in \nn$ and $b \in \nn$ satisfy $2 \le a < b$ and $\gcd(a, b) = 1$, set $q := a/b \in \qq$, and suppose that $b$ has at least two distinct prime divisors. Let $c \in \qq_{ > 0}$ satisfy $\den(c) \mid b^\infty$, and let $s \in \nn_0$ satisfy $a^s \mid \num(c)$. If $c > q^sB_q$ then the following statements hold.
	
	\begin{enumerate}
		\item[(1)] $c \in S_q$.
		\smallskip
		
		\item[(2)] For every $t \in S_q$ and every $\epsilon \in \rr_{ > 0}$, there are a nonunit $v \in T_q$ with $1 < v < 1 + \epsilon$ and an element $c' \in S_q$ such that $\den(c') \mid b^\infty$, $a^s \mid \num(c')$, $c' > q^sB_q$, and $t + c = v(t + c')$.
		\smallskip
		
		\item[(3)] For every $t \in S_q$, the element $t + c$ is not an atom of $S_q$.
		
	\end{enumerate}
\end{lemma}

\begin{proof}
	(1) Choose an integer $n > s$ sufficiently large that $\den(c) \mid b^n$. Since $a^s \mid \num(c)$ and $\gcd(a, b) = 1$, we can write $c = Ta^s/b^n$ with $T \in \nn$. We see that
	\[
	T = \frac{cb^n}{a^s} > b^{n - s} \cdot \frac{(a - 1)b}{b - a} > \frac{b^{n - s + 1}(a - 1) - a^{n - s + 1}(b - 1)}{b - a}.
	\]
	Hence the Frobenius number formula of Ong and Ponomarenko~\cite[Theorem]{OP08} implies that $T$ belongs to the numerical monoid $\langle a^{n - s}, a^{n - s - 1}b, \dots, b^{n - s} \rangle$. Thus,
	\[
	Ta^s \in \langle a^n, a^{n - 1}b, \dots, b^n \rangle,
	\]
	which implies that $c = Ta^s/b^n \in S_q$.
	\smallskip
	
	(2) Fix $t \in S_q$ and $\epsilon \in \rr_{ > 0}$. By Lemma~\ref{lem:tech-approximation}, for every $\eta \in \rr_{ > 0}$ there exists a nonunit $v \in T_q$ with $1 < v < 1 + \eta$ and $q^s \mid_{S_q} v - 1$.
	
	Choose such a $v = d_3/d_4 \in T_q$, where $d_3 \in \nn$ and $d_4 \in \nn$ are relatively prime, close enough to $1$ that $v < 1 + \epsilon$ and the rational number
	\[
	c' := \frac{t + c}{v} - t = (t + c)\frac{d_4}{d_3} - t
	\]
	satisfies $c' > q^sB_q$; this is possible because $c' \to c$ as $v \to 1$. Write $t + c = d_1/d_2$, where $d_1 \in \nn$ and $d_2 \in \nn$ are relatively prime. Since $t + c \in S_q$ and $v \in T_q$, the integers $d_2$, $d_3$, and $d_4$ have only prime divisors from $b$, and so they are relatively prime to $a$.
	
	We use congruences modulo $a^s$ after clearing denominators, which is legitimate because the relevant denominators are relatively prime to $a$. Since $a^s \mid \num(c)$ and $c = (t + c) - t$, we see that
	\[
	d_1 \equiv td_2 \pmod{a^s}.
	\]
	Since $q^s \mid_{S_q} v - 1$, we also have
	\[
	d_3 \equiv d_4 \pmod{a^s}.
	\]
	Therefore
	\[
	d_1d_4 - td_2d_3 \equiv 0 \pmod{a^s}.
	\]
	Now
	\[
	c' = \frac{d_1d_4}{d_2d_3} - t,
	\]
	and the denominator of $c'$ divides $b^\infty$. The preceding congruence gives $a^s \mid \num(c')$. By part~(1), applied to $c'$, we obtain $c' \in S_q$. Finally, $t + c = v(t + c')$, as required.
	\smallskip
	
	(3) Fix $t \in S_q$. Choose $\epsilon \in \rr_{ > 0}$ so small that no real number $v \in (1, 1 + \epsilon)$ equals $t + c$; if $t + c \le 1$ then every positive $\epsilon$ works. Part~(2) gives the factorization
	\[
	t + c = v(t + c')
	\]
	with $v \ne 1$ and $t + c' = (t + c)/v \ne 1$. In the present standing case, $S_q^*$ has no units other than $1$, so both factors are nonunits. Hence $t + c$ cannot be an atom.
\end{proof}

The preceding lemma yields a characterization of the atoms of $S_q^*$ when $q < 1/2$.


\begin{lemma} \label{lem:atoms-below-half}
	Let $q \in \qq \cap (0, 1/2)$ satisfy $\num(q) > 1$, and suppose that $\den(q)$ has at least two distinct prime divisors. Then the atoms of $S_q^*$ are precisely the elements with exactly one polynomial representation $f \in \nn_0[x]$, where $f$ is irreducible in $\nn_0[x]$.
\end{lemma}

\begin{proof}
	Set $a := \num(q) \in \nn$ and $b := \den(q) \in \nn$. Suppose that an element $c \in S_q^*$ has a unique polynomial representation $f \in \nn_0[x]$ and that $f$ is irreducible in $\nn_0[x]$. Then $c \ne 1$, since otherwise the constant polynomial $1$ would also represent $c$. If $c = c_1c_2$ with $c_1 \in S_q^*$ and $c_2 \in S_q^*$, choose representations $g_1 \in \nn_0[x]$ and $g_2 \in \nn_0[x]$ of these two factors. Uniqueness gives $f = g_1g_2$, so irreducibility forces $g_1 = 1$ or $g_2 = 1$. Thus, $c$ is an atom.
	
	We now prove that these are the only atoms. Suppose that there exists a polynomial representation $f \in \nn_0[x]$ such that the largest coefficient of $f$ is at least $a$ (by Lemma~\ref{lem:unique-polynomial-representation}, this is equivalent to there being more than one polynomial representation). We prove that $f(q)$ is not an atom.
	
	Choose $k \in \nn_0$ such that the coefficient of $x^k$ in $f$ is at least $a$. Fix a large integer $N \in \nn$, and consider the polynomial $g \in \zz[x]$ defined by
	\[
	g(x) := f(x) + (bx - a)(x^k + 2x^{k + 1} + 4x^{k + 2} + \dots + 2^Nx^{k + N}),
	\]
	which has nonnegative coefficients since $2a < b$. Choose $N$ large enough so that $k + N + 1 > \deg f$ and $2^Nb > B_q$. Then write
	\[
	f(q) = g(q) = t + 2^Nbq^{k + N + 1}
	\]
	for some $t \in S_q$. Set $s := k + N + 1 \in \nn$ and $c := 2^Nbq^s \in \qq_{ > 0}$. Since $\gcd(a, b) = 1$, reduction of this fraction gives $\den(c) \mid b^\infty$ and $a^s \mid \num(c)$. Moreover, $c > q^sB_q$. Lemma~\ref{lem:technical-division}(3) therefore shows that $f(q) = t + c$ is not an atom.
	
	Consequently, every atom has a unique polynomial representation $f$. This representation must be irreducible: if $f = gh$ with $g \in \nn_0[x] \setminus \{1\}$ and $h \in \nn_0[x] \setminus \{1\}$ then neither $g(q)$ nor $h(q)$ can equal $1$, since replacing that factor by the constant polynomial $1$ would give a different representation of $f(q)$. Hence $g(q)h(q)$ would be a factorization into nonunits, a contradiction.
\end{proof}

We now obtain the following corollary regarding elements with more than one polynomial representation.

\begin{corollary} \label{cor:small-technical-divisors}
	Let $q \in \qq \cap (0, 1/2)$ satisfy $\num(q) > 1$, with $\den(q)$ having at least two distinct prime divisors. For every $r \in S_q^*$ with more than one polynomial representation and every $\epsilon \in \rr_{ > 0}$, there is a nonunit $v \in T_q$ with $1 < v < 1 + \epsilon$ such that $v \mid_{S_q}r$.
\end{corollary}

\begin{proof}
	The construction in the proof of Lemma~\ref{lem:atoms-below-half} writes $r = t + c$, where $t \in S_q$ and $c \in \qq_{ > 0}$ satisfies the hypotheses of Lemma~\ref{lem:technical-division}. Part~(2) of Lemma~\ref{lem:technical-division} gives a nonunit divisor $v \in T_q$ with $1 < v < 1 + \epsilon$.
\end{proof}

The same construction provides a common technical divisor for finitely many elements with multiple representations.

\begin{lemma} \label{lem:finite-common-technical-divisor}
	Let $S_q$ be generated by $q \in \qq \cap (0, 1/2)$, where $\num(q) > 1$ and $\den(q)$ is not a prime power. Let $n \in \nn$, and let $r_1, \dots, r_n \in S_q^*$ each have more than one polynomial representation. For every $\epsilon \in \rr_{ > 0}$ there exists a nonunit $v \in T_q$ with $1 < v < 1 + \epsilon$ such that $v \mid_{S_q}r_i$ for every $i \in \ldb 1, n\rdb$.
\end{lemma}

\begin{proof}
	Fix $\epsilon \in \rr_{ > 0}$, and set $a := \num(q) \in \nn$ and $b := \den(q) \in \nn$. For each $i \in \ldb 1, n\rdb$, choose a polynomial representation $f_i \in \nn_0[x]$ of $r_i$ with some coefficient at least $a$; this is possible by Lemma~\ref{lem:unique-polynomial-representation}. Choose an index $k_i \in \nn_0$ such that the coefficient of $x^{k_i}$ in $f_i$ is at least $a$.
	
	Choose $N \in \nn$ large enough that $2^Nb > B_q$ and $k_i + N + 1 > \deg f_i$ for every $i \in \ldb 1, n\rdb$. For each $i \in \ldb 1, n \rdb$, define polynomials $h_i \in \zz[x]$ and $g_i \in \zz[x]$ by
	\[
	h_i(x) := x^{k_i} + 2x^{k_i + 1} + \cdots + 2^Nx^{k_i + N}
	\]
	and
	\[
	g_i(x) := f_i(x) + (bx - a)h_i(x).
	\]
	Since $q < 1/2$, we have $b > 2a$. Thus, the coefficient of $x^{k_i}$ in $g_i$ is still nonnegative, the intermediate new coefficients are $2^{j - 1}(b - 2a) > 0$ for $j \in \ldb 1, N \rdb$, and the final new coefficient is $2^Nb > 0$. Hence $g_i \in \nn_0[x]$ and $g_i(q) = f_i(q) = r_i$.
	
	For each $i \in \ldb 1, n \rdb$, put $s_i := k_i + N + 1 \in \nn$. Then
	\[
	r_i = g_i(q) = t_i + c_i,
	\]
	where $c_i := 2^Nbq^{s_i} \in \qq_{ > 0}$ and $t_i \in S_q$. The element $c_i$ satisfies $\den(c_i) \mid b^\infty$, $a^{s_i} \mid \num(c_i)$, and
	\[
	c_i = 2^Nbq^{s_i} > q^{s_i}B_q.
	\]
	Let $S := \max\{s_1, \dots, s_n\} \in \nn$.
	
	We now choose one technical divisor that works for all $i$. By Lemma~\ref{lem:tech-approximation}, for every $\eta \in \rr_{ > 0}$ there is a nonunit $v \in T_q$ such that $1 < v < 1 + \eta$ and $q^S \mid_{S_q}v - 1$. For each $i \in \ldb 1, n \rdb$, the rational number
	\[
	c_i'(v) := \frac{r_i}{v} - t_i
	\]
	tends to $c_i$ as $v \to 1$. Therefore after taking $\eta \in \rr_{ > 0}$ small enough, we may choose such a $v$ with $v < 1 + \epsilon$ and with $c_i'(v) > q^{s_i}B_q$ for every $i$.
	
	Since $S \ge s_i$, the divisibility $q^S \mid_{S_q}v - 1$ implies $q^{s_i} \mid_{S_q}v - 1$. The congruence argument in Lemma~\ref{lem:technical-division}(2), applied to $t_i + c_i$ and this common value of $v$, shows that each $c_i'(v)$ belongs to $S_q$. Consequently,
	\[
	r_i = v(t_i + c_i'(v))
	\]
	for every $i \in \ldb 1, n\rdb$. Thus, $v$ is a nonunit element of $T_q$ dividing every $r_i$, and $1 < v < 1 + \epsilon$.
\end{proof}

We conclude this subsection by showing that the characterization in Lemma~\ref{lem:atoms-below-half} does not extend to rational parameters $q \in \qq \cap (1/2, 1)$. The example works throughout this interval, including the non-ACCP cases considered above.

\begin{example}
	Let $q \in \qq \cap (1/2, 1)$, and set $a := \num(q) \in \nn$ and $b := \den(q) \in \nn$. Thus, $2 \le a < b < 2a$ and $\gcd(a, b) = 1$. The standing restriction on the number of prime divisors of $b$ is not needed for this example. We claim that the element $r \in S_q$ defined by
	\[
	r := 1 + (a - 1)q + aq^3
	\]
	is an atom of $S_q$ that has infinitely many polynomial representations.
	
	\noindent\textsc{Claim.} The only polynomial representations of $r$ are $1 + (a - 1)x + ax^3$ and
	\[
	1 + (a - 1)x + (b - a)(x^4 + x^5 + \dots + x^{n - 1}) + bx^n
	\]
	for integers $n \ge 4$.
	
	\noindent\textsc{Proof of Claim.} Let $f \in \nn_0[x]$ be a polynomial representation of $r$ different from $1 + (a - 1)x + ax^3$. By Gauss's lemma, there exists $c_0 + c_1x + \dots + c_dx^d \in \zz[x]$ (where $d \in \nn_0$, $c_0, \dots, c_d \in \zz$, and $c_d \neq 0$) such that
	\[
	f(x) = 1 + (a - 1)x + ax^3 + (bx - a)(c_0 + c_1x + \dots + c_dx^d).
	\]
	Extend the coefficient sequence by setting $c_i = 0$ for every integer $i > d$. If $c_i < 0$ for some $i \in \nn_0$ then nonnegativity of the coefficient of $x^{i + 1}$ in $f$ forces $c_{i + 1} < 0$: otherwise that coefficient is at most $(b - 1) - b < 0$, since every coefficient of $1 + (a - 1)x + ax^3$ is less than $b$. Iterating would contradict $c_{d + 1} = 0$. Hence every $c_i$ is nonnegative.
	
	Comparing coefficients of degrees $0$, $1$, $2$, and $3$ now gives $c_0 = c_1 = c_2 = 0$ and $c_3 \in \{0, 1\}$. For every $i \in \nn_0$ with $i \ge 3$, the coefficient of $x^{i + 1}$ in $f$ is $bc_i - ac_{i + 1} \ge 0$. Thus, $c_i = 0$ forces $c_{i + 1} = 0$, while $c_i = 1$ forces $c_{i + 1} \le b/a < 2$. Since the correction polynomial is nonzero, its nonzero coefficients form a consecutive block $c_3 = \cdots = c_d = 1$. This gives the claimed expression with $n = d + 1$. Conversely, each listed polynomial has nonnegative coefficients and evaluates to $r$, by $bq = a$. For $n = 4$, the sum $x^4 + \cdots + x^{n - 1}$ is empty.
	\hfill $\square$
	
	It remains to prove that every listed representation is irreducible in $\nn_0[x]$. Every nontrivial factorization has factors with constant coefficient $1$, hence both factors have positive degree. The degree-$3$ representation cannot factor: the degrees would be $1$ and $2$, but the positive leading coefficient of the quadratic factor would contribute to its zero $x^2$ coefficient. The degree-$4$ representation cannot factor either. Its nonzero linear coefficient forces a linear term in one factor; its zero coefficients in degrees $2$ and $3$ force the other nonconstant factor to have degree at least $4$, giving total degree at least $5$.
	
	Now let $n \in \nn$ satisfy $n \ge 5$, and define $f \in \nn_0[x]$ by
	\[
	f(x) := 1 + (a - 1)x + (b - a)(x^4 + \cdots + x^{n - 1}) + bx^n
	\]
	and suppose $f = gh$ with nonunits $g \in \nn_0[x]$ and $h \in \nn_0[x]$. The zero coefficient of $x^2$ prevents both factors from having positive linear coefficients. After interchanging them if needed, coefficient comparison gives
	\begin{align*}
		[x^0]g & = [x^0]h = 1, & [x^1]g & = a - 1, & [x^1]h & = 0, \\
		[x^2]g & = [x^3]g = 0, & [x^2]h & = [x^3]h = 0.
	\end{align*}
	Put $d := \deg g \in \nn$ and $e := \deg h \in \nn$, so $e \ge 4$ and $d + e = n$.
	
	If $d > 1$ then $d \ge 4$ and $e + 1 < n$. Consequently,
	\[
	a - 1 \le (a - 1)[x^e]h \le [x^{e + 1}]f = b - a \le a - 1.
	\]
	All inequalities are equalities, so $[x^e]h = 1$. Comparing leading coefficients gives $[x^d]g = b$. But $4 \le d < n$, and the constant coefficient of $h$ is $1$, so
	\[
	b - a = [x^d]f \ge [x^d]g = b,
	\]
	a contradiction.
	
	Thus, $d = 1$ and $g = 1 + (a - 1)x$. Write $h = 1 + x^4m(x)$ with $m \in \nn_0[x]$. Substitution into $f = gh$ and comparison of the terms of degree at least $4$ give the following identity for the polynomial $k \in \nn_0[x]$:
	\[
	k(x) := (b - a)(1 + x + \cdots + x^{n - 5}) + bx^{n - 4} = (1 + (a - 1)x)m(x).
	\]
	If $n \ge 6$, comparison of constant and linear coefficients gives
	\[
	[x^0]m = b - a, \qquad
	b - a = (a - 1)(b - a) + [x^1]m.
	\]
	Since $b - a > 0$ and $a \ge 2$, this forces $a = 2$ and $[x^1]m = 0$. However,
	\[
	k(-1) =
	\begin{cases}
		b, & n\text{ even}, \\
		-a, & n\text{ odd},
	\end{cases}
	\]
	contradicting divisibility by $1 + x$. If $n = 5$ then $m$ is constant, so $(a - 1)(b - a) = b$. The equality $\gcd(b - a, b) = 1$ forces $b - a = 1$, and the same equation then yields $b = a - 1$, again a contradiction.
	
	All representations are therefore irreducible. Since $S_q^*$ has only the unit $1$, every factorization of $r$ into two nonunits would give a reducible polynomial representation of $r$. Thus, $r$ is an atom. The listed representations have distinct degrees, so there are infinitely many of them. \hfill $\blacksquare$
\end{example}

\medskip
\subsection{The Furstenberg Property}

We now construct rational monogenic semidomains that fail the Furstenberg property. These examples are non-atomic, since every atomic monoid is Furstenberg.

By Theorem~\ref{thm:ff-bf-accp}, $S_q^*$ satisfies the FF property in each of the following cases. It is therefore atomic and, in particular, Furstenberg in each listed case.
\begin{enumerate}
	\item[(1)] $q \in \nn^{-1}$.
	\smallskip
	
	\item[(2)] $q \ge 1$.
	\smallskip
	
	\item[(3)] $q \in (0, 1)$ such that $\num(q) > 1$ and $\den(q)$ is a prime power.
\end{enumerate}
In the remaining results on non-Furstenberg examples, we restrict our attention to the case where $q \in (0, 1)$ such that $\num(q) > 1$ and $\den(q)$ is not a prime power.

We proceed by reducing the problem of the Furstenberg property from $S_q^*$ to $T_q$.

\begin{theorem} \label{thm:furstenberg-reduction}
	Let $M$ be a commutative cancellative monoid, let $\omega \colon M \to (\nn_0, +)$ be a homomorphism, and set $M_0 := \omega^{-1}(0) \subseteq M$. Then $M$ is Furstenberg if and only if $M_0$ is Furstenberg. In particular, for every $q \in \qq_{ > 0}$, the monoid $S_q^*$ is Furstenberg if and only if $T_q$ is Furstenberg.
\end{theorem}

\begin{proof}
	The submonoid $M_0$ is divisor-closed because $\omega$ takes nonnegative values. Hence its units and atoms are exactly the units and atoms of $M$ that lie in $M_0$. It follows that the Furstenberg property of $M$ implies that of $M_0$.
	
	Conversely, suppose $M_0$ is Furstenberg, and assume that a nonunit $c \in M$ has no atom divisor. Put $k := \omega(c) \in \nn_0$. Starting with $r_0 := c \in M$, choose elements $z_i \in M$ and $r_{i + 1} \in M$ successively so that
	\[
	r_i = z_i r_{i + 1},
	\]
	for every $i \in \ldb 0, k - 1\rdb$, where both factors are nonunits, and put $z_k := r_k \in M$. Each step is possible because the residual nonunit $r_i$ divides $c$ and therefore is not an atom. We obtain $c = z_0 \cdots z_k$ as a product of $k + 1$ nonunits, with the same conclusion for $k = 0$ by taking $z_0 = c$. Since
	\[
	k = \omega(c) = \sum_{i = 0}^k \omega(z_i),
	\]
	there exists $j \in \ldb 0, k \rdb$ such that $z_j$ lies in $M_0$. The Furstenberg property of $M_0$ gives an atom of $M_0$ dividing $z_j$. This is an atom of $M$ dividing $c$, a contradiction. Apply the result to $M = S_q^*$ and $\omega = \omega_q$ for the final assertion.
\end{proof}

We now use the approximation of $1$ by technical fractions to give a sufficient condition for failure of the Furstenberg property.

\begin{prop} \label{prop:technical-atom-gap}
	Let $q \in \qq \cap (0, 1)$ satisfy $\num(q) > 1$, and suppose that $\den(q)$ has at least two distinct prime divisors. If $\inf\mathcal A(T_q) > 1$ then $S_q^*$ is not Furstenberg. In particular, this holds if $T_q$ has finitely many atoms.
\end{prop}

\begin{proof}
	The hypotheses violate condition~(e) of Theorem~\ref{thm:ff-bf-accp}, and $q$ is not a unit fraction. Condition~(d) therefore implies that $1$ is a limit point of $T_q$. Also, $T_q \subseteq \qq_{ \ge 1}$ and its only unit is $1$.
	
	If $\mathcal A(T_q)$ is nonempty, let $\lambda := \inf\mathcal A(T_q) \in \rr_{ > 1}$ and choose $r \in T_q$ with $1 < r < \lambda$. If $\mathcal A(T_q)$ is empty, choose an element $r \in T_q$ with $1 < r < 2$. In either case $r$ is a nonunit with no atom divisor in $T_q$: in the nonempty case, a factorization $r = uw$ with $u \in \mathcal A(T_q)$ and $w \in T_q$ would give $r \ge u \ge \lambda$. Divisor-closedness implies that every atom of $S_q^*$ dividing $r$ would be an atom of $T_q$. Thus, $r$ has no atom divisor in $S_q^*$.
\end{proof}

Proposition~\ref{prop:technical-atom-gap} yields explicit parameters for which $S_q^*$ is not Furstenberg.

\begin{theorem} \label{thm:non-furstenberg}
	Let $S_q$ be the monogenic semidomain generated by $q \in \qq_{ > 0}$. If $\den(q)$ has at least two distinct prime divisors and $1 < \num(q) < \min \{p \in \pp : p \mid \den(q) \}$ then $S_q$ does not satisfy the Furstenberg property.
\end{theorem}

\begin{proof}
	We prove that $T_q$ is antimatter in this case, which suffices by Proposition~\ref{prop:technical-atom-gap}. Set $a := \num(q) \in \nn$ and $b := \den(q) \in \nn$. Note that $b > a^2 \ge 2a$, so $q < 1/2$ and Lemma~\ref{lem:atoms-below-half} applies.
	
	Suppose, by way of contradiction, that $c \in \mathcal A(T_q)$. By divisor-closedness, $c \in \mathcal A(S_q^*)$, so Lemma~\ref{lem:atoms-below-half} gives a unique polynomial representation $f \in \nn_0[x]$. Moreover, $c > 1$ because $T_q \subseteq \qq_{ \ge 1}$ and $c$ is a nonunit. Write
	\[
	f(x) = c_0 + c_1x + c_2x^2 + \dots + c_kx^k,
	\]
	where $k := \deg f \in \nn_0$, $c_i \in \nn_0$ for every $i \in \ldb 0, k \rdb$, $c_k > 0$, and $c_i < a$ for every such $i$, by Lemma~\ref{lem:unique-polynomial-representation}. Then $c = N/b^k$, where $N := \sum_{i = 0}^k c_i a^i b^{k - i} \in \nn$. Choose a prime $p \in \pp$ dividing $\num(c)$; it exists because $c > 1$. Since $c \in T_q$, we have $p \mid b$, and $p \mid N$ because the reduced numerator divides $N$. Reducing modulo $p$ gives $p \mid c_ka^k$. As $\gcd(a, b) = 1$, this implies $p \mid c_k$, which contradicts $0 < c_k < a < p$. This also covers $k = 0$, when $N = c_0$. Thus, $T_q$ has no atoms.
\end{proof}

The preceding theorem also gives the following failure of atomicity.

\begin{corollary}
	Let $S_q$ be the monogenic semidomain generated by $q \in \qq_{ > 0}$. If $\den(q)$ has at least two distinct prime divisors and $1 < \num(q) < \min \{p \in \pp : p \mid \den(q) \}$ then $S_q$ is not atomic.
\end{corollary}

\begin{proof}
	By Theorem~\ref{thm:non-furstenberg}, $S_q^*$ is not Furstenberg. Since every atomic monoid is Furstenberg, $S_q^*$ is not atomic.
\end{proof}

The antimatter conclusion in the proof of Theorem~\ref{thm:non-furstenberg} admits the following extension, in which the prime divisors of $\den(q)$ are replaced by the full prime-power factors of $\den(q)$. The hypothesis is the one that will guarantee MCD-finiteness in Theorem~\ref{thm:mcd-finite-prime-powers}.

\begin{theorem} \label{thm:finite-technical-atoms}
	Let $a \in \nn$ and $b \in \nn$ satisfy $2 \le a < b$ and $\gcd(a, b) = 1$, and set $q := a/b \in \qq$. Suppose that $b$ has at least two distinct prime divisors and that $p^{\nu_p(b)} > a$ for every prime $p \in \pp$ dividing $b$. Then the following statements hold.
	
	\begin{enumerate}
		\item[(1)] For $c \in T_q$, the element $c$ has a unique polynomial representation if and only if $c \in \ldb 1, a - 1 \rdb$.
		\smallskip
		
		\item[(2)] $\mathcal A(T_q) = \{p \in \pp : p \mid b \text{ and } p < a\}$. In particular, $T_q$ has only finitely many atoms.
		\smallskip
		
		\item[(3)] $S_q^*$ is not Furstenberg, and so it is not atomic.
		
	\end{enumerate}
\end{theorem}

\begin{proof}
	Since $b$ is the product of at least two prime-power factors exceeding $a$, we have $b > a^2 \ge 2a$. Thus, $q < 1/2$, and Lemma~\ref{lem:atoms-below-half} applies.
	
	(1) Let $c \in T_q$ have exactly one polynomial representation $f = \sum_{i = 0}^k c_i x^i \in \nn_0[x]$, where $k := \deg f \in \nn_0$ and $c_i \in \nn_0$ for every $i \in \ldb 0, k \rdb$. By Lemma~\ref{lem:unique-polynomial-representation}, $0 \le c_i < a$ for every $i$, and $c_k \ge 1$. Set $N := \sum_{i = 0}^k c_i a^i b^{k - i} \in \nn$, so that $c = N/b^k$. Since $\num(c) \mid b^\infty$ and $N = \num(c) \cdot b^k/\den(c)$, every prime divisor of $N$ divides $b$. Fix a prime $p \in \pp$ dividing $b$, and put $e_p := \nu_p(b) \in \nn$. Every term of $N$ other than $c_k a^k$ is divisible by $b$, whence $N \equiv c_k a^k \pmod{p^{e_p}}$. Since $\gcd(a, p) = 1$ and $1 \le c_k < a < p^{e_p}$, we have $\nu_p(c_k a^k) = \nu_p(c_k) < e_p$, and therefore $\nu_p(N) = \nu_p(c_k)$. Because $N$ is supported on the prime divisors of $b$,
	\[
	N = \prod_{p \mid b} p^{\nu_p(N)} = \prod_{p \mid b} p^{\nu_p(c_k)},
	\]
	which divides $c_k$. Hence $N \le c_k$. On the other hand, $N \ge c_k a^k$, so $k = 0$ and $c = N = c_0 \in \ldb 1, a - 1 \rdb$. Conversely, every integer in $\ldb 1, a - 1 \rdb$ that belongs to $T_q$ has a unique polynomial representation, namely the constant polynomial, by Lemma~\ref{lem:unique-polynomial-representation}.
	\smallskip
	
	(2) By divisor-closedness, $\mathcal A(T_q) = T_q \cap \mathcal A(S_q^*)$. Let $c \in \mathcal A(T_q)$. By Lemma~\ref{lem:atoms-below-half}, $c$ has a unique polynomial representation, which is irreducible in $\nn_0[x]$. By part~(1), this representation is the constant polynomial $c \in \ldb 1, a - 1 \rdb$. A constant polynomial is irreducible in $\nn_0[x]$ precisely when it is a prime number, so $c$ is a prime, and $c \mid b$ because $c \in T_q$. Conversely, let $p \in \pp$ divide $b$ and satisfy $p < a$. Then $p \in T_q$, the constant polynomial $p$ is the unique representation of $p$ by Lemma~\ref{lem:unique-polynomial-representation}, and it is irreducible in $\nn_0[x]$. Thus, $p$ is an atom of $S_q^*$ by Lemma~\ref{lem:atoms-below-half}, and hence an atom of $T_q$.
	\smallskip
	
	(3) By part~(2), $\mathcal A(T_q)$ is finite, so $S_q^*$ is not Furstenberg by Proposition~\ref{prop:technical-atom-gap}. Since every atomic monoid is Furstenberg, $S_q^*$ is not atomic.
\end{proof}

The following examples illustrate the preceding results and the possible presence of technical atoms outside their hypotheses.

\begin{example} \label{ex:technical-atoms}
	The semidomains $S_{2/15}$, $S_{2/35}$, and $S_{3/55}$ are non-Furstenberg. Therefore they are non-atomic rational monogenic semidomains. Theorem~\ref{thm:finite-technical-atoms} applies to $q = 5/72$ as well: here $\mathcal A(T_q) = \{2, 3\}$, and $S_{5/72}^*$ is not Furstenberg, although Theorem~\ref{thm:non-atomic} does not apply because $5 > 3$.
	
	When the hypothesis of Theorem~\ref{thm:non-furstenberg} is not satisfied, the technical monoid can contain atoms, as Table~\ref{tab:technical-atoms} shows. Both parameters in the table also violate the hypothesis of Theorem~\ref{thm:finite-technical-atoms}, since $4 < 5$ for $q = 5/12$ and $8 < 11$ for $q = 11/24$; this is consistent with the presence of technical atoms of positive degree. All displayed coefficients are less than $\num(q)$, so the representations are unique. The four linear polynomials are primitive and therefore irreducible in $\nn_0[x]$. The cubic $6x^3 + 4x^2 + 9x + 3$ is primitive and irreducible over $\qq$: modulo $7$ its values at $0, 1, \dots, 6$ are respectively $3$, $1$, $1$, $4$, $4$, $2$, and $6$, so its reduction is a cubic without a root. Lemma~\ref{lem:atoms-below-half} therefore makes all five values atoms of $S_q^*$. Their reduced numerators and denominators involve only primes dividing $\den(q)$, so they lie in $T_q$ and are atoms there as well. These examples establish the existence of technical atoms, not the Furstenberg property for the corresponding monoids.
	\begin{table}[ht]
		\centering \caption{Explicit atoms in technical monoids.} \label{tab:technical-atoms}
		
		\begin{tabular}{c|c|c}
			$q$ & Value & Polynomial representation \\ \hline $5/12$ & $9/4$ & $1 + 3q$ \\ & $8/3$ & $1 + 4q$\\ \hline $11/24$ & $27/8$ & $2 + 3q$\\ & $32/3$ & $7 + 8q$\\ & $2187/256$ & $3 + 9q + 4q^2 + 6q^3$
			
		\end{tabular}
		
	\end{table}
	\hfill $\blacksquare$
\end{example}

\medskip
\subsection{Atomicity}

We now construct rational monogenic semidomains that fail atomicity. These examples are broader than the non-Furstenberg rational monogenic semidomains constructed in Theorem~\ref{thm:non-furstenberg}, but they are not comparable with those of Theorem~\ref{thm:finite-technical-atoms}.

\begin{theorem} \label{thm:non-atomic}
	Let $S_q$ be the monogenic semidomain generated by $q \in \qq_{ > 0}$. If $\den(q)$ has at least two distinct prime divisors and $1 < \num(q) < \max \{p \in \pp : p \mid \den(q) \}$ then $S_q$ is not atomic.
\end{theorem}

\begin{proof}
	Suppose, by way of contradiction, that $S_q$ is atomic. Set $a := \num(q) \in \nn$ and $b := \den(q) \in \nn$. Let $p_0 := \max \{ p \in \pp : p \mid \den(q) \} \in \pp$. Note that $2a < 2p_0 \le b$, since $b$ has at least two distinct prime divisors; hence $q < 1/2$.
	
	Note that $p_0 \in S_q$. Because $S_q$ is atomic, the element $p_0$ admits a factorization $p_0 = a_1a_2 \dotsb a_k$, where $k \in \nn$ and $a_1, \dots, a_k \in \mathcal A(S_q^*)$. For each $i \in \ldb 1, k \rdb$, let $f_i(x) \in \nn_0[x]$ be a polynomial representation of $a_i$, and let $t_i \in \nn$ be the leading coefficient of $f_i$.
	
	Let $g(x) := f_1(x)f_2(x) \dotsb f_k(x) \in \nn_0[x]$, which is a polynomial representation of $p_0$. First, suppose that $g$ is constant. Then $p_0 = g(0) = f_1(0)f_2(0) \dotsb f_k(0)$, so there exists an index $m \in \ldb 1, k \rdb$ such that $p_0 \mid f_m(0)$. Hence the constant coefficient of $f_m$ is greater than $a$, because $p_0 > a$. Therefore $f_m(x) + bx - a \in \nn_0[x]$ is another polynomial representation of $a_m$, different from $f_m$. Lemma~\ref{lem:atoms-below-half} then implies that $a_m$ is not an atom, a contradiction.
	
	On the other hand, suppose $g$ is nonconstant. We see that $g(q) = p_0$, and by Gauss's lemma we obtain
	\[
	g(x) = p_0 + (bx - a)h(x)
	\]
	for some polynomial $h(x) \in \zz[x]$. Let $c \in \zz \setminus \{0\}$ be the leading coefficient of $h$. The leading coefficient of $g$ is both $bc$ and $t_1t_2 \dotsb t_k$. In particular, $p_0 \mid bc = t_1t_2 \dotsb t_k$. Therefore there exists an index $n \in \ldb 1, k \rdb$ such that $p_0 \mid t_n$. It follows that $t_n \ge p_0 > a$, so $f_n(x) + (bx - a)x^{\deg f_n} \in \nn_0[x]$ is another polynomial representation of $a_n$, different from $f_n$. Lemma~\ref{lem:atoms-below-half} then implies that $a_n$ is not an atom, a contradiction.
\end{proof}

\medskip
\subsection{Strong Atoms}

We conclude this section by determining the strong atoms of atomic rational monogenic semidomains. Recall that an atom $u \in \mathcal A(S_q^*)$ is a \emph{strong atom} if, for every $n \in \nn$, the power $u^n$ has a unique factorization up to order and associates. For atoms and strong atoms in the additive monoids of related semidomains, see Dani et al.~\cite{DDGLPVZ26}.

We first deal with the factorial cases. If $q \in \nn \cup \nn^{-1}$ then $S_q^*$ is factorial by Deng, Gotti, and Zeng~\cite{DGZ26}. Hence every atom of $S_q^*$ is strong. Moreover, $S_q^*$ has infinitely many nonassociate atoms: if $q \in \nn$ then $S_q = \nn_0$, whose atoms are the ordinary primes, while if $q = 1/m$ for some $m \in \nn$, every prime not dividing $m$ is an atom of $S_q^*$. Therefore throughout the remainder of this subsection, we write $q = a/b$ in lowest terms, where $a \in \nn$ and $b \in \nn$ satisfy $a > 1$ and $b > 1$.

The next two lemmas establish that $q$ is strong and that every other atom eventually admits each technical element as a divisor of one of its powers.

\begin{lemma} \label{lem:q-strong-atom}
	Let $q = a/b \in \qq_{ > 0}$ be written in lowest terms, where $a \in \nn$ and $b \in \nn$ satisfy $a > 1$ and $b > 1$. Then $q$ is a strong atom of $S_q^*$.
\end{lemma}

\begin{proof}
	For every $n \in \nn$, the polynomial $x^n$ is the unique representation of $q^n$, by Lemma~\ref{lem:unique-polynomial-representation}. Its only factorization into atoms of $\nn_0[x]^*$ is the product of $n$ copies of $x$. Proposition~\ref{prop:unique-divisors} therefore shows that $q^n$ has the unique factorization consisting of $n$ copies of $q$.
\end{proof}

\begin{lemma} \label{lem:strong-eventual-b}
	Let $q = a/b \in \qq_{ > 0}$ be in lowest terms, with $a \in \nn$ and $b \in \nn$ satisfying $a > 1$ and $b > 1$. For every $u \in \mathcal A(S_q^*) \setminus \{q\}$ and every $t \in T_q$, the divisibility $t \mid_{S_q^*} u^k$ holds for all sufficiently large $k \in \nn$. In particular, this holds for $t = b$.
\end{lemma}

\begin{proof}
	Choose a representation $f \in \nn_0[x]$ of $u$, and put $d := \deg f \in \nn_0$, $A := b^d f(q) \in \nn$, and $M := \max\{a, b\} \in \nn$. If $d > 0$ then $f(0) > 0$: otherwise $q$ divides the atom $u$, forcing $u = q$ because the monoid is reduced. Thus,
	\[
	A \ge a^d + b^d > M^d.
	\]
	If $d = 0$ then $A = u > 1 = M^d$, so the same strict inequality holds.
	
	Write $t = N/D$, where $N \in \nn$ and $D \in \nn$ are relatively prime, and choose $\ell \in \nn$ with $N \mid b^\ell$. For $j \in \nn$, define the numerical monoid $M_j := \langle a^j, a^{j - 1}b, \dots, b^j\rangle$. Formula~\eqref{eq:geometric-frobenius} gives a constant $C \in \rr_{ > 0}$, depending only on $a$ and $b$, such that $F(M_j) \le C M^{j + 1}$. Since $A/M^d > 1$, for all sufficiently large $k \in \nn$,
	\[
	A^k > C M^{dk + \ell + 1} \ge F(M_{dk + \ell}).
	\]
	Consequently, $A^k \in M_{dk + \ell}$, so $u^k/b^\ell \in S_q$. Therefore
	\[
	\frac{u^k}{t} = D\cdot \frac{b^\ell}{N}\cdot \frac{u^k}{b^\ell} \in S_q^*.
	\]
	This proves the assertion.
\end{proof}

The two lemmas yield the following classification when $S_q^*$ is atomic.

\begin{theorem} \label{thm:strong-atoms}
	Suppose $q = a/b \in \qq_{ > 0}$ is written in lowest terms, where $a \in \nn$ and $b \in \nn$ satisfy $a > 1$ and $b > 1$. Then the following statements hold.
	
	\begin{enumerate}
		\item[(1)] If $b = p^e$ for some prime $p \in \pp$ and $e \in \nn$ then $\mathcal A_{\mathrm{str}}(S_q^*) = \{q, p\}$.
		\smallskip
		
		\item[(2)] If $b$ has at least two distinct prime divisors and $S_q^*$ is atomic then $\mathcal A_{\mathrm{str}}(S_q^*) = \{q\}$.
		
	\end{enumerate}
\end{theorem}

\begin{proof}
	By Lemma~\ref{lem:q-strong-atom}, $q$ is a strong atom of $S_q^*$.
	
	(1) Suppose $b = p^e$ for some prime $p \in \pp$ and $e \in \nn$. Recall from Theorem~\ref{thm:ff-bf-accp} that $S_q^*$ satisfies the FF property and is therefore atomic.
	
	We first show that $p$ is strong. Recall from Section~\ref{sec:algebra} that $T_q$ is divisor-closed in $S_q^*$ and that $T_q \subseteq \qq_{ \ge 1}$ because $a > 1$ (see also Deng, Gotti, and Zeng~\cite[Proposition~4.9(1)]{DGZ26}). As every element of $T_q$ is supported on the single prime $p$, we obtain $T_q = \{1, p, p^2, \dots\}$. The prime $p$ is the unique atom of $T_q$. For every $n \in \nn$, divisor-closedness implies that every atom of $S_q^*$ dividing $p^n$ lies in $T_q$. Hence $p^n$ factors uniquely as a product of $n$ copies of $p$.
	
	It remains to show no other strong atoms exist. Let $u \in \mathcal A_{\mathrm{str}}(S_q^*) \setminus \{q\}$. By Lemma~\ref{lem:strong-eventual-b}, we have $u^k = bv$ for a sufficiently large $k \in \nn$ and some $v \in S_q^*$.
	
	Since $S_q^*$ is atomic, factor $b$ and $v$ into atoms, yielding a factorization of $u^k$. Because $u$ is strong, all atoms appearing must be associate to $u$. Since $a > 1$ and $b > 1$, $S_q^*$ is reduced, so $b = u^m$ for some $m \in \nn$. Writing $u = h/d$, where $h \in \nn$ and $d \in \nn$ are relatively prime, the equality $b = u^m$ forces $d^m \mid h^m$, so $d = 1$ and $u \in \nn$. If $u$ were composite in $\nn$, its integer factors greater than $1$ would be nonunits of the reduced monoid $S_q^*$, contradicting that $u$ is an atom. Thus, $u$ is an ordinary prime. The only prime divisor of $b$ is $p$, so $u = p$ as desired.
	\smallskip
	
	(2) Under the assumption that $S_q^*$ is atomic, we may proceed as above, using Lemma~\ref{lem:strong-eventual-b} to obtain $b = u^m$ for every strong atom $u \neq q$, with $m \in \nn$. As above, $b = u^m$ forces $u \in \nn$, and irreducibility of $u$ then forces $u$ to be an ordinary prime, which is a contradiction since $b$ is not a prime power.
\end{proof}

The preceding theorem leaves open the strong atoms in the non-atomic cases. We conclude this section with the following problem.

\begin{prob}
	Determine the set of strong atoms of $S_q$, where $q = a/b \in \qq_{ > 0}$ is in lowest terms, $a \in \nn$ and $b \in \nn$ satisfy $a > 1$ and $b > 1$, and $b$ is not a prime power.
\end{prob}

\bigskip
\section{GCDs and MCDs}
\label{sec:gcd-mcd}

In this section, we study divisibility properties of the multiplicative monoid $S_q^*$ of a rational monogenic semidomain. We first determine the parameters $q$ for which $S_q^*$ satisfies the GCD and quasi-GCD properties. We then turn to the weaker MCD and MCD-finite properties, for which several cases have been solved, although the general problem remains open. Following the notation of Deng, Gotti, and Zeng~\cite{DGZ26}, throughout this section let $S_q$ denote the monogenic semidomain generated by $q \in \qq_{ > 0}$. Recall that the technical monoid is the divisor-closed submonoid
\[
T_q = \{c \in S_q^* : \num(c) \mid \den(q)^\infty\}.
\]

\medskip
\subsection{The GCD and Quasi-GCD Properties}

We begin with the stronger GCD and quasi-GCD properties. Recall that the quasi-GCD property states that every nonempty finite subset of $S_q^*$ has at most one MCD up to associates.

\begin{theorem} \label{thm:gcd-classification}
	Let $S_q$ be the monogenic semidomain generated by a positive rational parameter $q$. The following conditions are equivalent.
	
	\begin{enumerate}
		\item[(a)] $S_q$ is GCD.
		\smallskip
		
		\item[(b)] $S_q$ is quasi-GCD.
		\smallskip
		
		\item[(c)] $q \in \nn \cup \nn^{-1}$.
	\end{enumerate}
\end{theorem}

\begin{proof}
	(a) $\Rightarrow$ (b): In a GCD monoid, every MCD of a finite nonempty set is associate to its GCD.
	\smallskip
	
	(b) $\Rightarrow$ (c): Suppose that $q = a/b$ is in lowest terms, where $a$ and $b$ are integers greater than $1$. The polynomials $1 + x$, $1 + x^3$, $1 + x + x^2$, and $1 + x^2 + x^4$ are irreducible in $\nn_0[x]$: each has constant coefficient $1$ and a prime value at $1$, so a factor with value $1$ must be the constant polynomial $1$. Their coefficients are less than $\min\{a, b\}$, and Lemma~\ref{lem:atom-condition-sq} shows that their evaluations at $q$ are atoms of $S_q^*$.
	
	Consider the finite subset of $S_q^*$ given by
	\[
	E := \{1 + q + q^2 + q^3 + q^4 + q^5, \ 1 + q + q^3 + q^4\}.
	\]
	The identities
	\[
	1 + q + q^2 + q^3 + q^4 + q^5 = (1 + q^3)(1 + q + q^2) = (1 + q)(1 + q^2 + q^4)
	\]
	and $1 + q + q^3 + q^4 = (1 + q)(1 + q^3)$ show that $u := 1 + q \in S_q^*$ and $v := 1 + q^3 \in S_q^*$ are common divisors of $E$. The quotient set $E/u$ consists of the two nonassociate atoms $1 + q^2 + q^4$ and $1 + q^3$; they are distinct because their difference is $q^2(q^2 - q + 1) > 0$. Likewise, $E/v$ consists of the distinct atoms $1 + q + q^2$ and $1 + q$. Thus, $u$ and $v$ are MCDs. Since $q \ne 1$, we have $u \ne v$, and reducedness makes them nonassociate. Hence $S_q^*$ is not quasi-GCD.
	\smallskip
	
	(c) $\Rightarrow$ (a): If the generating parameter $q$ belongs to $\nn \cup \nn^{-1}$ then $S_q^*$ is factorial by Corollary~\ref{cor:uf-ap}, and it is well known that every factorial monoid has the GCD property.
\end{proof}

The following example makes the distinction between greatest and maximal common divisors explicit.

\begin{example}
	Let $q := 2/3 \in \qq$, and set $r := 665/243 \in S_q^*$ and $s := 175/81 \in S_q^*$. The identities
	\[
	r = 1 + q + q^2 + q^3 + q^4 + q^5 \qquad \text{and} \qquad s = 1 + q + q^3 + q^4
	\]
	give unique polynomial representations, since every coefficient is less than $2$. Set $u := 1 + q = 5/3 \in S_q^*$ and $v := 1 + q^3 = 35/27 \in S_q^*$. Then
	\[
	\{r/u, s/u\} = \{133/81, 35/27\} \qquad \text{and} \qquad \{r/v, s/v\} = \{19/9, 5/3\}.
	\]
	The first quotient set consists of the evaluations of $1 + x^2 + x^4$ and $1 + x^3$, and the second consists of the evaluations of $1 + x + x^2$ and $1 + x$. Each polynomial has constant coefficient $1$ and prime value at $1$, so it is irreducible in $\nn_0[x]$. Lemma~\ref{lem:atom-condition-sq} makes both quotient sets pairs of distinct atoms. Thus, $u$ and $v$ are MCDs of $\{r, s\}$. They are distinct and nonassociate because $S_{2/3}^*$ is reduced. In particular, $\{r, s\}$ has no GCD. \hfill $\blacksquare$
\end{example}

\medskip
\subsection{The MCD and MCD-Finite Properties}

We now study the MCD and MCD-finite properties. The preceding factorization results settle the ACCP cases. The remaining parameters lie in $(0, 1)$, have numerator greater than $1$, and have denominators divisible by at least two distinct primes. We prove MCD-finiteness and failure of the MCD property under additional hypotheses on these denominator primes and their powers.

\begin{prop}
	If $q \in [1, + \infty) \cup \nn^{-1}$ or $\den(q)$ is a prime power then $S_q^*$ is both MCD and MCD-finite.
\end{prop}

\begin{proof}
	Recall that the ACCP implies the MCD property. Hence by Theorem~\ref{thm:ff-bf-accp}, $S_q^*$ satisfies the MCD property for all of the stated parameters.
	
	Moreover, Theorem~\ref{thm:ff-bf-accp} shows that $S_q^*$ satisfies the FF property. Therefore every element of $S_q^*$ has only finitely many divisors up to associates. Let $R \subseteq S_q^*$ be finite and nonempty, and fix $r \in R$. Every common divisor of $R$ is a divisor of $r$, so $R$ has only finitely many common divisors up to associates, and hence only finitely many MCDs up to associates. Thus, $S_q^*$ is MCD-finite.
\end{proof}

Recall the digit-numerator set $\mathcal D_{a, b}$ from Section~\ref{sec:background}. We first show that, under the following hypotheses, this set omits infinitely many primes.

\begin{lemma} \label{lem:prime-exist}
	Let $q = a/b \in (0, 1) \setminus \nn^{-1}$ be a rational number in lowest terms, with $a$ and $b$ in $\nn$, such that $b$ has at least two distinct prime divisors and all prime divisors of $b$ are greater than $a$. Then infinitely many primes do not lie in $\mathcal{D}_{a, b}$.
\end{lemma}

\begin{proof}
	Let $X \in \rr$ with $X \ge 1$. We prove that the number of integers in $[1, X]$ belonging to $\mathcal{D}_{a, b}$ is $O(X^{1 - \epsilon}\log X)$ for a fixed real number $\epsilon > 0$.
	
	Choose $A \in \mathcal{D}_{a, b}$ with $A \le X$, and write $A = c_0b^n + \cdots + c_na^n$, where $n \in \nn_0$, $c_i \in \ldb 0, a - 1 \rdb$ for every $i \in \ldb 0, n \rdb$, and $c_n > 0$. Let $j \in \ldb 0, n \rdb$ be the least index with $c_j \ne 0$, and set $k := n - j \in \nn_0$. Then the term $c_ja^jb^{n - j}$ is present in $A$, so
	\[
	\left(\dfrac{b}{a}\right)^k \le a^{n - k}b^k \le X.
	\]
	
	Since $c_n > 0$, we have $n \le \log_a X = O(\log X)$. For fixed $n$ and $k$, each coefficient $c_i$ with $i \in \ldb 0, n - k - 1 \rdb$ vanishes, while there are at most $a^{k + 1}$ choices for the remaining coefficients. Since $k \le \log_{b/a} X$, summing over the possible indices gives the following bound for the number of elements in $[1, X] \cap \mathcal{D}_{a, b}$:
	\[
	O(\log X) \sum_{0 \le k \le \log_{b/a} X} a^{k + 1}.
	\]
	Since $a$ is fixed, the summation on the right is $O\left(a^{\log_{b/a} X}\right) = O(X^{\log_{b/a}a}).$
	
	Since $b$ has at least two distinct prime divisors, each greater than $a$, we have $b > a^2$ and $\log_{b/a} a < 1$. Set $\epsilon := 1 - \log_{b/a} a \in \rr_{ > 0}$. The preceding bound is therefore
	\[
	O(X^{1 - \epsilon}\log X) = o(X/\log X).
	\]
	Since the number of primes in $[1, X]$ is asymptotic to $X/\log X$ by the prime number theorem, we conclude that infinitely many primes do not lie in $\mathcal{D}_{a, b}$.
\end{proof}

The next lemma excludes unique polynomial representations for the quotients of these primes by technical divisors.

\begin{lemma} \label{lem:irreducible}
	Let $q = a/b \in (0, 1) \setminus \nn^{-1}$ be a rational number in lowest terms, with $a$ and $b$ in $\nn$, such that $b$ has at least two distinct prime divisors and all prime divisors of $b$ are greater than $a$. If $P \in \pp \setminus \mathcal{D}_{a, b}$ does not divide $b$ then for each $d \in T_q$ dividing $P$, the element $P/d$ does not have a unique polynomial representation in $S_q$.
\end{lemma}

\begin{proof}
	Suppose by way of contradiction that $P/d$ had a unique polynomial representation. Since $b > a^2$ and $a > 1$, we have $q < 1/2$. By Lemma~\ref{lem:unique-polynomial-representation} we can write
	\[
	\frac{P}{d} = \dfrac{c_0b^n + c_1ab^{n - 1} + \dots + c_na^n}{b^n}
	\]
	for some $n \in \nn_0$, with $c_i \in \ldb 0, a - 1 \rdb$ for every $i \in \ldb 0, n - 1 \rdb$ and $c_n \in \ldb 1, a - 1 \rdb$. Denote the numerator by $A \in \nn$. By definition, $A \in \mathcal{D}_{a, b}$, and $A \equiv c_na^n \pmod{b}$. Since $\gcd(a, b) = 1$ and all prime divisors of $b$ are greater than $a > c_n$, it follows that $A$ is coprime to $b$. Now
	\[
	d = \dfrac{P\cdot b^n}{A} \in T_q,
	\]
	and since $\nu_P(d) = 0$, we conclude $P \mid A$. Write $A = A_0P$ with $A_0 \in \nn$. Then $d = b^n/A_0, $ where $\gcd(A_0, b) = \gcd(A, b) = 1.$ Since the denominator of $d$ divides $b^\infty$, we obtain $A_0 = 1$ and $A = P$. This contradicts $P \notin \mathcal{D}_{a, b}$.
\end{proof}

We combine these results to construct a set with no MCD.

\begin{theorem} \label{thm:large-prime-non-mcd}
	Let $q = a/b \in (0, 1) \setminus \nn^{-1}$ be a rational number in lowest terms, with $a$ and $b$ in $\nn$, such that $b$ has at least two distinct prime divisors and all prime divisors of $b$ are greater than $a$. Then $S_q$ does not satisfy the MCD property.
\end{theorem}

\begin{proof}
	From Lemma~\ref{lem:prime-exist}, choose distinct sufficiently large primes $R$ and $S$ not contained in $\mathcal{D}_{a, b}$ and not dividing $b$.
	
	Suppose that $d \in S_q^*$ is a common divisor of $\{R, S\}$. We show that $d$ cannot be maximal. If a prime $p \in \pp$ with $p \nmid b$ divided $\num(d)$ then $d \mid_{S_q}R$ would imply $p \mid R$, while $d \mid_{S_q}S$ would imply $p \mid S$, because reduced denominators of elements of $S_q$ involve only prime divisors of $b$. Since $R$ and $S$ are distinct primes, this is impossible. Hence every prime divisor of $\num(d)$ divides $b$, and so $d \in T_q$.
	
	From Lemma~\ref{lem:irreducible}, we know that $R/d$ and $S/d$ have more than one polynomial representation. Applying Lemma~\ref{lem:finite-common-technical-divisor}, there are nonunit elements of $T_q$ arbitrarily close to $1$ dividing both $R/d$ and $S/d$, implying that $d$ is not an MCD. Thus, the set $\{R, S\}$ has no MCD.
\end{proof}

The following example gives an explicit pair with no MCD.

\begin{example}
	Set $q := 2/15 \in \qq$ and $E := \{7, 11\} \subseteq S_q^*$. Every member of $\mathcal D_{2, 15}$ below $17$ is a power of $2$. Indeed, let $A := \sum_{i = 0}^n c_i2^i15^{n - i} \in \nn$, where $n \in \nn_0$, $c_i \in \{0, 1\}$ for every $i \in \ldb 0, n \rdb$, and $c_n = 1$. If $c_j = 1$ for some $j \in \ldb 0, n - 1 \rdb$ then
	\[
	A \ge 2^n + 2^j15^{n - j} \ge 2^n + 15\cdot 2^{n - 1} \ge 17.
	\]
	Otherwise, $A = 2^n$. Consequently, $7$ and $11$ do not belong to $\mathcal D_{2, 15}$, and neither prime divides $15$. Every common divisor $d \in S_q^*$ of $E$ is technical. Lemma~\ref{lem:irreducible} shows that both $7/d$ and $11/d$ have more than one representation, so Lemma~\ref{lem:finite-common-technical-divisor} supplies a further nonunit common divisor. Thus, $E$ has no MCD. Nevertheless, $S_q^*$ is MCD-finite by Proposition~\ref{prop:mcd-finite-cat}; the set of MCDs of $E$ is empty, hence finite. \hfill $\blacksquare$
\end{example}

The next lemma locates a quotient with a unique polynomial representation. We then bound the possible uniquely represented divisors under a weaker hypothesis on the denominator.

\begin{lemma} \label{lem:mcd-finite-reduced}
	Let $q = a/b \in (0, 1) \setminus \nn^{-1}$ be a rational number in lowest terms, with $a$ and $b$ in $\nn$, such that $b$ has at least two distinct prime divisors and all prime divisors of $b$ are greater than $a$. If $d \in S_q^*$ is an MCD of a finite nonempty subset $R \subset S_q^*$ then there exists an element $r \in R$ such that $r/d$ has a unique polynomial representation.
\end{lemma}

\begin{proof}
	The argument is the same as in Theorem~\ref{thm:large-prime-non-mcd}. Suppose by way of contradiction that no element of $R/d$ has a unique polynomial representation. By Lemma~\ref{lem:finite-common-technical-divisor}, we can choose a nonunit element of $T_q$ dividing every element contained in $R/d$. This contradicts the maximality of $d$.
\end{proof}

To bound the possible quotients, we first prove the finite-divisor statement under the weaker hypothesis on the full prime-power factors of the denominator.

\begin{lemma} \label{lem:finite-divisor-prime-powers}
	Let $a$ and $b$ be relatively prime positive integers with $1 < a < b$, and set $q := a/b \in \qq$. Suppose that $p^{\nu_p(b)} > a$ for every prime $p$ dividing $b$. Then, for each $r \in S_q^*$, there are only finitely many $u \in S_q^*$ with exactly one polynomial representation such that $u \mid_{S_q^*} r$.
\end{lemma}

\begin{proof}
	Fix $r \in S_q^*$, and write $r = N/D$ in lowest terms, where $N$ and $D$ are positive integers. The denominator $D$ divides a power of $b$. Let $u \in S_q^*$ divide $r$ and have exactly one polynomial representation. If that representation is constant then $u$ is an integer in $\ldb 1, a - 1 \rdb$ by Lemma~\ref{lem:unique-polynomial-representation}, so there are only finitely many possibilities. Hence assume that the unique representation has degree $n \in \nn$. By the same lemma, we can write
	\[
	u = \frac{c_0b^n + c_1ab^{n - 1} + \cdots + c_n a^n}{b^n} = \frac{A}{b^n},
	\]
	where $A \in \nn$, $c_i \in \ldb 0, a - 1 \rdb$ for every $i \in \ldb 0, n - 1 \rdb$, and $c_n \in \ldb 1, a - 1 \rdb$. Fix a prime $p \in \pp$ dividing $b$, and set $e_p := \nu_p(b) \in \nn$. Reducing modulo $p^{e_p}$ gives
	\[
	A \equiv c_na^n \pmod{p^{e_p}},
	\]
	since every other term in $A$ is divisible by $b$. Since $\gcd(a, p) = 1$ and $1 \le c_n < a < p^{e_p}$, we have $\nu_p(c_na^n) = \nu_p(c_n) < e_p$. Therefore
	\[
	\nu_p(A) = \nu_p(c_n)
	\]
	for every prime $p$ dividing $b$.
	
	Set $g := \gcd(A, b^n) \in \nn$. Since $n \ge 1$, the valuation calculation gives
	\[
	g = \prod_{p \mid b} p^{\nu_p(c_n)}.
	\]
	In particular, $g \mid c_n$, and $A' := A/g \in \nn$ is relatively prime to $b$. Thus, the reduced expression for $u$ is
	\[
	u = \frac{A'}{b^n/g}.
	\]
	Since $u \mid_{S_q^*} r$, we have
	\[
	\frac{r}{u} = \frac{N(b^n/g)}{DA'} \in S_q^*.
	\]
	Every element of $S_q$ has reduced denominator dividing a power of $b$. Since $A'$ is relatively prime to $b$, its prime-power factors must cancel with $N$, and hence $A' \mid N$. Moreover, $g \mid c_n$ and $A \ge c_na^n$, so
	\[
	A' = A/g \ge a^n(c_n/g) \ge a^n.
	\]
	Consequently, $a^n \le A' \le N$, which bounds $n$ by a number depending only on $r$. For each fixed $n$, there are only finitely many choices for $c_0, c_1, \dots, c_n$, and these coefficients determine $u$. This proves the claim.
\end{proof}

The denominator-prime version is an immediate special case.

\begin{lemma} \label{lem:finite-divisor}
	Let $q = a/b \in (0, 1) \setminus \nn^{-1}$ be a rational number in lowest terms, with $a$ and $b$ in $\nn$, such that $b$ has at least two distinct prime divisors and all prime divisors of $b$ are greater than $a$. Then, for each $r \in S_q^*$, there are only finitely many $u \in S_q^*$ with exactly one polynomial representation such that $u \mid_{S_q^*} r$.
\end{lemma}

\begin{proof}
	Every prime $p$ dividing $b$ satisfies $p^{\nu_p(b)} \ge p > a$. Therefore Lemma~\ref{lem:finite-divisor-prime-powers} applies.
\end{proof}

The general finite-divisor criterion yields the following MCD-finiteness theorem.

\begin{theorem} \label{thm:mcd-finite-prime-powers}
	Let $q \in \qq \cap (0, 1)$ with $\num(q) > 1$ satisfy $p^{\nu_p(\den(q))} > \num(q)$ for all primes $p$ dividing $\den(q)$. Then the monogenic semidomain $S_q$ is MCD-finite.
\end{theorem}

\begin{proof}
	Write $q = a/b$ in lowest terms, where $a$ and $b$ are positive integers and $a > 1$. If $b$ is a prime power then $S_q^*$ is an FF monoid by Theorem~\ref{thm:ff-bf-accp}. Since $S_q^*$ is reduced in this case, each element has only finitely many divisors, and so every finite nonempty subset has only finitely many common divisors. Hence $S_q$ is MCD-finite. Thus, we may assume that $b$ has at least two distinct prime divisors.
	
	Under the present hypothesis, we have
	\[
	b = \prod_{p \mid b} p^{\nu_p(b)} > a^2 \ge 2a.
	\]
	Therefore $q < 1/2$, so Lemma~\ref{lem:finite-common-technical-divisor} applies. Let $R \subseteq S_q^*$ be finite and nonempty, and let $d \in S_q^*$ be an MCD of $R$. If every element of $R/d$ had more than one polynomial representation then that lemma would provide a nonunit common divisor of $R/d$, contradicting the maximality of $d$. Hence there exists $r \in R$ such that $r/d$ has exactly one polynomial representation, as in the proof of Lemma~\ref{lem:mcd-finite-reduced}.
	
	For each fixed $r \in R$, the element $r/d$ divides $r$ in $S_q^*$. Lemma~\ref{lem:finite-divisor-prime-powers} therefore gives only finitely many possibilities for $r/d$, and hence for $d = r/(r/d)$. Since $R$ is finite, it has only finitely many MCDs. Thus, $S_q$ is MCD-finite.
\end{proof}

We recover the criterion involving the prime divisors themselves as a special case.

\begin{prop} \label{prop:mcd-finite-cat}
	Let $q = a/b \in (0, 1) \setminus \nn^{-1}$ be a rational number in lowest terms, with $a$ and $b$ in $\nn$, such that $b$ has at least two distinct prime divisors, and all prime divisors of $b$ are greater than $a$. Then $S_q$ is MCD-finite.
\end{prop}

\begin{proof}
	Every prime $p$ dividing $b$ satisfies $p^{\nu_p(b)} \ge p > a$. The conclusion therefore follows from Theorem~\ref{thm:mcd-finite-prime-powers}.
\end{proof}

Combining the nonexistence and finiteness results gives the following distinction.

\begin{corollary}
	Let $q = a/b \in (0, 1) \cap \qq$ be in lowest terms, where $a$ and $b$ are positive integers and $a > 1$. Suppose that $b$ has at least two distinct prime divisors, all greater than $a$. Then $S_q^*$ is MCD-finite but is not MCD.
\end{corollary}

\begin{proof}
	This follows by combining Theorem~\ref{thm:large-prime-non-mcd} and Proposition~\ref{prop:mcd-finite-cat}.
\end{proof}

MCD-finiteness does not impose a uniform bound on the number of MCDs. The following construction applies to every nonfactorial rational monogenic semidomain.

\begin{prop} \label{thm:unbounded-mcd-cardinality}
	Let $q \in \qq_{ > 0} \setminus (\nn \cup \nn^{-1})$. For every $n \in \nn$, there exist $r$ and $s$ in $S_q^*$ such that $\{r, s\}$ has at least $n$ pairwise nonassociate MCDs.
\end{prop}

\begin{proof}
	Fix $n \in \nn$. Choose distinct odd primes $m_1, \dots, m_n \in \pp$ such that $m_i > m_1 + \cdots + m_{i - 1}$ for every $i \in \ldb 2, n \rdb$, and set $L := 2m_1 \cdots m_n \in \nn$. Define the polynomials $A$ and $B$ in $\nn_0[x]$ by
	\[
	A(x) := 1 + x + \cdots + x^{L - 1}
	\]
	and
	\[
	B(x) := \prod_{j = 1}^n (1 + x^{m_j}).
	\]
	The choice of the $m_i$ makes their subset sums distinct, so both $A$ and $B$ have coefficients at most $1$. Write $q = a/b$ in lowest terms, where $a$ and $b$ are integers greater than $1$. Lemma~\ref{lem:unique-polynomial-representation} therefore gives unique representations for $A(q)$ and $B(q)$.
	
	Fix $i \in \ldb 1, n \rdb$, and define the polynomial $D_i \in \nn_0[x]$ by $D_i(x) := 1 + x^{m_i}$. This polynomial divides $B$ in $\nn_0[x]$ and also divides $A$, since
	\[
	\frac{A(x)}{D_i(x)} = (1 + x + \cdots + x^{m_i - 1})\sum_{j = 0}^{L/(2m_i) - 1}x^{2m_i j} \in \nn_0[x].
	\]
	We claim that $D_i$ is an MCD of $\{A, B\}$ in $\nn_0[x]^*$. The polynomial $A = (x^L - 1)/(x - 1)$ is squarefree over $\qq$, whereas
	\[
	B(x) = (1 + x)^n\prod_{j = 1}^n\Phi_{2m_j}(x).
	\]
	Consequently, the monic GCD in $\zz[x]$ is
	\[
	\gcd_{\zz[x]}(A, B) = (1 + x)\prod_{j = 1}^n\Phi_{2m_j}(x).
	\]
	If a common divisor $E \in \nn_0[x]$ were a proper multiple of $D_i$ then $E$ would be monic, because it divides the monic polynomial $A$. Unique factorization in $\zz[x]$ would force $E/D_i$ to be a product of a nonempty collection of the factors $\Phi_{2m_j}$ with $j \ne i$. Each such factor is $1 - x + x^2 - \cdots + x^{m_j - 1}$, so the coefficient of $x$ in that product is the negative of the number of factors. This contradicts $E/D_i \in \nn_0[x]$. Hence $D_i$ is an MCD.
	
	Proposition~\ref{prop:unique-divisors} transfers this maximality to $D_i(q)$ in $S_q^*$. Indeed, a common divisor of $A(q)$ and $B(q)$ has a unique polynomial representation, which divides both $A$ and $B$, and the same correspondence preserves divisibility by $D_i$. Finally, the numbers $D_i(q)$ are distinct because $q > 0$, $q \ne 1$, and the $m_i$ are distinct. Since $S_q^*$ is reduced, these MCDs are pairwise nonassociate.
\end{proof}

The next example lies beyond the hypothesis that every denominator prime exceeds the numerator.

\begin{example}
	The monogenic semidomain $S_{3/20}$ satisfies the MCD-finite property by Theorem~\ref{thm:mcd-finite-prime-powers}, since the full prime-power factors of $20$ are $4$ and $5$, both greater than $3$. \hfill $\blacksquare$
\end{example}

We conclude this section with a reduction of the MCD problem to technical common divisors.

For a finite nonempty subset $R \subseteq S_q^*$, recall the common-divisor set $\CD_{S_q^*}(R)$, and set
\[
\CD_{T_q}(R)
:= \CD_{S_q^*}(R) \cap T_q.
\]

\begin{prop}
	Let $q \in \qq_{ > 0}$. Then the monogenic semidomain $S_q$ satisfies the MCD property if and only if the following condition holds: for every finite nonempty subset $R \subseteq S_q^*$, if $d_0 \in \CD_{S_q^*}(R)$ has maximal $\omega_q$ value then $\CD_{T_q}(R/d_0)$ has a maximal element under divisibility.
\end{prop}

\begin{proof}
	First suppose that the stated condition holds. Fix a finite nonempty subset $R \subseteq S_q^*$. Such an element $d_0$ exists because $1$ is a common divisor and $\omega_q$ takes nonnegative integer values bounded above by $\omega_q(r)$ for a fixed $r \in R$. Suppose that $u \in S_q^* \setminus T_q$ divided $R/d_0$ in $S_q$. Then $d_0u$ would be a common divisor of $R$ with $\omega_q(d_0u) > \omega_q(d_0)$, contradicting the choice of $d_0$. Therefore
	\[
	\CD_{S_q^*}(R/d_0) = \CD_{T_q}(R/d_0).
	\]
	By assumption, this set has a maximal element $t \in T_q$. Suppose $d_0t \mid_{S_q^*} c$, where $c \in S_q^*$ is another common divisor of $R$ not associate to $d_0t$. Then
	\[
	\dfrac{c}{d_0} = t \cdot \dfrac{c}{td_0} \in S_q
	\]
	is a larger divisor of $R/d_0, $ contradicting the maximality of $t$. It follows that $d_0t$ is an MCD of $R$.
	
	Conversely, suppose that $S_q$ satisfies the MCD property, and let $R \subseteq S_q^*$ be finite and nonempty. Choose a common divisor $d_0 \in S_q^*$ with maximal $\omega_q$ value. The same argument above yields $\CD_{S_q^*}(R/d_0) = \CD_{T_q}(R/d_0)$. As a consequence, an MCD of $R/d_0$ in $S_q$ must remain a maximal element.
\end{proof}

\bigskip
\section{ACCP and Quasi-ACCP}
\label{sec:accp}

Theorem~\ref{thm:ff-bf-accp} classifies the rational monogenic semidomains satisfying the ACCP. We now study the weaker quasi-ACCP property and characterize the individual elements satisfying the ACCP when the entire monoid fails it.

\medskip
\subsection{The Quasi-ACCP Property}

We first iterate the technical-division lemma to construct nonstabilizing chains of principal ideals.

\begin{lemma} \label{lem:large-coefficient-non-accp}
	Consider the monogenic semidomain $S_q$ generated by $q \in (0, 1) \cap \qq$, and assume that $\num(q) > 1$ and $\den(q)$ is not a prime power. Let $c \in \qq_{ > 0}$ satisfy $\den(c) \mid \den(q)^\infty$, and suppose $s \in \nn_0$ satisfies $\num(q)^s \mid \num(c)$. If $c > q^s B_q$ then for every $t \in S_q$ the element $t + c \in S_q$ does not satisfy the ACCP.
\end{lemma}

\begin{proof}
	Fix $t \in S_q$, and set $c_0 := c \in \qq_{ > 0}$. By Lemma~\ref{lem:technical-division}(2), we can recursively choose nonunits $v_i \in T_q$ and elements $c_i \in S_q$ for $i \in \nn$ satisfying the same denominator, numerator-divisibility, and lower-bound hypotheses as $c_0$ such that
	\[
	t + c_{i - 1} = v_i(t + c_i)
	\]
	for every $i \in \nn$. Thus, we obtain the ascending chain
	\[
	(t + c_0)S_q^* \subseteq (t + c_1)S_q^* \subseteq (t + c_2)S_q^* \subseteq \cdots.
	\]
	Each containment is strict. Indeed, equality at the $i$th step would give $t + c_i = w(t + c_{i - 1})$ for some $w \in S_q^*$. Combining the two identities and canceling the nonzero element $t + c_i$ gives $v_iw = 1$, contradicting the choice of $v_i$. Hence this chain does not stabilize, and $t + c$ does not satisfy the ACCP.
\end{proof}

The next two examples illustrate the coefficient condition for failure of the ACCP.

\begin{example}
	Let $q := 3/10 \in \qq$. Then $\num(q) = 3$, $\den(q) = 10$ is not a prime power, and
	\[
	B_q = \frac{\num(q) - 1}{1 - q} = \frac{20}{7}.
	\]
	Take $s := 0 \in \nn_0$ and $c := 4 \in \qq_{ > 0}$. Then $\den(c) = 1 \mid 10^\infty$, $3^0 \mid \num(c)$, and $4 > 20/7$. Hence Lemma~\ref{lem:large-coefficient-non-accp}, with $t = 0$, shows that the element $4 \in S_{3/10}$ does not satisfy the ACCP. Notice also that $4 \in T_q$, since the only prime dividing its numerator is $2$, which divides $\den(q)$. \hfill $\blacksquare$
\end{example}

\begin{example}
	Let $q := 2/15 \in \qq$. Then $\num(q) = 2$, $\den(q) = 15$ is not a prime power, and
	\[
	B_q = \frac{\num(q) - 1}{1 - q} = \frac{15}{13}.
	\]
	Take $s := 1 \in \nn_0$ and $c := 15q = 2 \in \qq_{ > 0}$. We have $\den(c) = 1 \mid 15^\infty$, $2^1 \mid \num(c)$, and
	\[
	c = 2 > qB_q = \frac{2}{13}.
	\]
	Therefore Lemma~\ref{lem:large-coefficient-non-accp} shows, for example with $t = 1$, that
	\[
	1 + 15q = 3
	\]
	does not satisfy the ACCP in $S_{2/15}$. This illustrates that failure of the ACCP can result from a large coefficient of a low power of $q$, even when the value of the element is small. \hfill $\blacksquare$
\end{example}

For parameters below $1/2$, unique polynomial representations also characterize the ACCP elements.

\begin{lemma} \label{lem:accp-unique-representation}
	Let $q \in \qq \cap (0, 1/2)$ such that $\num(q) > 1$ and $\den(q)$ is not a prime power. Then the elements of $S_q^*$ satisfying the ACCP are precisely the elements that have exactly one polynomial representation.
\end{lemma}

\begin{proof}
	If an element has more than one polynomial representation, the proof of Lemma~\ref{lem:atoms-below-half} writes it as $t + c$, where $t \in S_q$, $c \in \qq_{ > 0}$, and some $s \in \nn_0$ satisfy the hypotheses of Lemma~\ref{lem:large-coefficient-non-accp}. The latter lemma shows that it does not satisfy the ACCP. Conversely, every element with a unique polynomial representation satisfies the ACCP by Proposition~\ref{prop:unique-divisors}.
\end{proof}

We next use this characterization to prove that the quasi-ACCP and MCD properties coincide for these parameters.

\begin{prop} \label{prop:quasi-accp-iff-mcd}
	Let $q \in \qq \cap (0, 1/2)$ such that $\num(q) > 1$ and $\den(q)$ is not a prime power. Then $S_q$ is quasi-ACCP if and only if $S_q$ is MCD.
\end{prop}

\begin{proof}
	The quasi-ACCP property implies the MCD property for every commutative cancellative monoid, as recalled in Section~\ref{sec:background}. For the converse, assume that $S_q^*$ is MCD, and let $E \subseteq S_q^*$ be finite and nonempty. Choose an MCD $d \in S_q^*$ of $E$. Since $S_q^*$ is reduced, the quotient set $E/d$ has no common divisor other than $1$.
	
	If every element of $E/d$ had more than one polynomial representation, Lemma~\ref{lem:finite-common-technical-divisor} would supply a nonunit technical common divisor. Hence some element of $E/d$ has a unique polynomial representation. By Lemma~\ref{lem:accp-unique-representation}, this element satisfies the ACCP. Thus, $d$ and this element verify the quasi-ACCP condition for $E$.
\end{proof}

The following conjecture would rule out the quasi-ACCP and MCD properties for every remaining parameter below $1/2$.

\begin{conjecture} \label{conj:technical-atom-gap}
	Let $q \in \qq \cap (0, 1/2)$ such that $\num(q) > 1$ and $\den(q)$ is not a prime power. Then $\inf(\mathcal{A}(T_q)) > 1$, where $\inf\emptyset := + \infty$.
\end{conjecture}

We first prove a lemma independent of this conjecture and then derive its consequence.

\begin{lemma} \label{lem:small-not-furst}
	Let $q \in \qq \cap (0, 1/2)$ such that $\num(q) > 1$ and $\den(q)$ is not a prime power. If there exist two elements $a \in S_q^*$ and $b \in S_q^*$ such that $a \nmid_{S_q} b$, $b \nmid_{S_q} a$, and neither $a$ nor $b$ has an atom divisor then $S_q$ is neither MCD nor quasi-ACCP.
\end{lemma}

\begin{proof}
	Set $S := \{a, b\} \subseteq S_q^*$. We claim that there is no common divisor $d \in S_q^*$ of $S$ for which some element of $S/d$ satisfies the ACCP. This shows that $S_q$ does not satisfy the quasi-ACCP and, by Proposition~\ref{prop:quasi-accp-iff-mcd}, is not MCD.
	
	Assume by way of contradiction that there exists $d \in S_q^*$ such that $d$ is a common divisor of $S$ and some element of $S/d$ satisfies the ACCP. Without loss of generality, suppose $a/d$ satisfies the ACCP. Note that $d \ne a$ due to the condition $a \nmid_{S_q} b$, so $a/d \ne 1$. Then $a/d$ factors into atoms, so $a/d$ is divisible by at least one atom. This contradicts our original assumption, so $S_q$ is not quasi-ACCP, as desired.
\end{proof}

The conjectured gap above $1$ supplies the two incomparable elements required by the preceding lemma.

\begin{prop} \label{prop:conj-implies-non-q-accp}
	Let $q \in \qq \cap (0, 1/2)$ such that $\num(q) > 1$ and $\den(q)$ is not a prime power. If Conjecture~\ref{conj:technical-atom-gap} holds for this $q$ then the semidomain $S_q$ is neither MCD nor quasi-ACCP.
\end{prop}

\begin{proof}
	Write $q = m/n$ in lowest terms, where $m \in \nn$ and $n \in \nn$. Choose a real number $c$ with $1 < c < \inf\mathcal A(T_q)$. Since $T_q \subseteq \qq_{ \ge 1}$ and is divisor-closed, no element of $T_q \cap (1, c)$ has an atom divisor in $S_q^*$. By Lemma~\ref{lem:small-not-furst}, it therefore suffices to find two divisibility-incomparable elements in this interval.
	
	Write $n = rs$ with relatively prime integers $r \in \nn$ and $s \in \nn$ satisfying $r \ge 2$ and $s \ge 2$, and choose $\ell \in \nn$ such that $q^\ell B_q < c - 1$. The numbers $\log r$ and $\log s$ are linearly independent over $\qq$: an integer relation would give $r^u = s^v$ for nonzero integers $u$ and $v$, contrary to coprimality. The irrational-rotation density theorem then shows that $\{u\log r - v\log s : u, v \in \nn\}$ is dense in $\rr$. For completeness, positive $u$ can make $u\log r$ modulo $\log s$ arbitrarily close to a prescribed residue, with $u$ arbitrarily large; subtracting the appropriate positive multiple $v\log s$ gives an approximation to every real target.
	
	Choose $u \in \nn$ and $v \in \nn$ so that the rational number
	\[
	a := \left(\frac{r^u}{s^v}\right)^{\phi(m^\ell)}
	\]
	satisfies $1 + q^\ell B_q < a < c$. Both $\num(a)$ and $\den(a)$ divide $n^\infty$, and Euler's theorem gives $a \equiv 1 \pmod{m^\ell}$. Here and below, rational congruences are interpreted by inverting denominators, all of which are relatively prime to $m$. Thus, $m^\ell \mid \num(a - 1)$ and $a - 1 > q^\ell B_q$. Lemma~\ref{lem:technical-division}(1), applied to $a - 1$, gives $a - 1 \in S_q$. Hence $a \in T_q \cap (1, c)$.
	
	Since $a \ne 1$, choose $k \in \nn$ with $k > \ell$ such that $a \not \equiv 1 \pmod{m^k}$. By the same density argument, choose $x \in \nn$ and $y \in \nn$ so that the rational number
	\[
	b := \left(\frac{r^x}{s^y}\right)^{\phi(m^k)}
	\]
	satisfies $a < b < \min\{c, a + aq^k\}$. Again, the numerator and denominator of $b$ are supported on the prime divisors of $n$, and $b \equiv 1 \pmod{m^k}$. Since $k > \ell$ and $b - 1 > q^\ell B_q$, applying Lemma~\ref{lem:technical-division}(1) to $b - 1$ with exponent $\ell$ gives $b \in T_q \cap (1, c)$.
	
	If $b$ divided $a$ in $S_q^*$ then additivity of $\omega_q$ would give $a/b \in T_q$, contradicting $a/b < 1$. Therefore $b \nmid_{S_q^*} a$. Suppose, conversely, that $a$ divided $b$. Then $b/a \in T_q$, so every polynomial representation has a positive constant coefficient. Since
	\[
	1 < \frac{b}{a} < 1 + q^k < 2,
	\]
	the constant coefficient must be $1$, and every other nonzero coefficient must have exponent greater than $k$. It follows that $b/a \equiv 1 \pmod{m^k}$. However,
	\[
	\frac{b}{a} \equiv a^{-1} \not \equiv 1 \pmod{m^k},
	\]
	a contradiction. Thus, $a \nmid_{S_q^*} b$, and Lemma~\ref{lem:small-not-furst} completes the proof.
\end{proof}

The finite-atom classification yields an unconditional consequence for the prime-power range.

\begin{corollary} \label{cor:prime-power-non-mcd}
	Let $q = a/b \in \qq \cap (0, 1)$ be in lowest terms, where $a \in \nn$ and $b \in \nn$, with $a > 1$ and $b$ divisible by at least two distinct primes. If $p^{\nu_p(b)} > a$ for every prime divisor $p$ of $b$ then $S_q^*$ is MCD-finite but is neither MCD nor quasi-ACCP.
\end{corollary}

\begin{proof}
	The hypothesis gives $b > a^2$, so $q < 1/2$. Theorem~\ref{thm:finite-technical-atoms} identifies $\mathcal A(T_q)$ with the finite set of primes dividing $b$ that are less than $a$. Its infimum is therefore greater than $1$, including the possibility of an empty atom set. Proposition~\ref{prop:conj-implies-non-q-accp} now rules out the MCD and quasi-ACCP properties. The MCD-finite property follows from Theorem~\ref{thm:mcd-finite-prime-powers}.
\end{proof}

This consequence contains the earlier result based on large denominator primes.

\begin{remark}
	Under the hypotheses of Corollary~\ref{cor:prime-power-non-mcd}, if every prime divisor of $b$ is greater than $a$ then the atom set of $T_q$ is empty. Thus, the preceding corollary also gives another proof of Theorem~\ref{thm:large-prime-non-mcd}.
\end{remark}

\medskip
\subsection{A Coefficient Criterion for ACCP Elements}

We now remove the restriction $q < 1/2$ and characterize ACCP elements throughout the non-ACCP parameter range.

\begin{theorem} \label{thm:accp-coefficient-bound}
	Let $q \in (0, 1) \cap \qq$ satisfy $\num(q) > 1$, and suppose that $\den(q)$ has at least two distinct prime divisors. For an element $r \in S_q^*$, the following conditions are equivalent.
	
	\begin{enumerate}
		\item[(a)] The element $r$ satisfies the ACCP.
		\smallskip
		
		\item[(b)] There is a constant bounding every coefficient of every polynomial representation of $r$.
		\smallskip
		
		\item[(c)] Every coefficient of every polynomial representation of $r$ is at most $B_q = (\num(q) - 1)/(1 - q)$.
	\end{enumerate}
\end{theorem}

\begin{proof}
	(a) $\Rightarrow$ (c): If a representation has a coefficient $C \in \nn$ with $C > B_q$ at $x^s$ for some $s \in \nn_0$, write $r = t + Cq^s$ with $t \in S_q$. Reduction of $Cq^s$ leaves its denominator supported on the prime divisors of $\den(q)$ and its numerator divisible by $\num(q)^s$. Thus, Lemma~\ref{lem:large-coefficient-non-accp} shows that $r$ does not satisfy the ACCP, a contradiction.
	\smallskip
	
	(c) $\Rightarrow$ (b): The constant $B_q$ is a uniform bound.
	\smallskip
	
	(b) $\Rightarrow$ (a): Write $a := \num(q) \in \nn$ and $b := \den(q) \in \nn$. For each $z \in S_q^*$ whose representation coefficients are bounded, let $m(z) \in \nn$ be the largest coefficient occurring among all its polynomial representations. This maximum exists because the coefficients are nonnegative integers. We first establish two claims.
	
	\noindent\textsc{Claim 1.} Let $u \in S_q^*$ and $v \in S_q^*$. If $u \mid_{S_q^*} v$ and $m(v)$ exists then $m(u)$ exists and $m(u) \le m(v)$.
	
	\noindent\textsc{Proof of Claim 1.} Fix a representation $g \in \nn_0[x]$ of the nonzero quotient $v/u$, and choose a positive coefficient of $g$. Multiplying a representation $f \in \nn_0[x]$ of $u$ by $g$ gives a representation of $v$ with a coefficient at least as large as each coefficient of $f$. Thus, every coefficient of every representation of $u$ is at most $m(v)$. \hfill $\blacksquare$
	
	\noindent\textsc{Claim 2.} Let $u \in S_q^*$ and $v \in S_q^*$. Suppose $u \mid_{S_q^*} v$, $m(v)$ exists, $m(u) \ge a$, and $v/u$ is not a nonnegative integral power of $q$. Then $m(u) < m(v)$.
	
	\noindent\textsc{Proof of Claim 2.} Choose a representation $f \in \nn_0[x]$ of $u$ with a coefficient equal to $m(u)$, and choose a representation $g \in \nn_0[x]$ of $v/u$. If a coefficient of $g$ is at least $2$ then the product $fg$ has a coefficient at least $2m(u)$, which proves the claim. Otherwise, every coefficient of $g$ is at most $1$, and $g$ has at least two nonzero terms because the quotient is not a power of $q$. Thus, $fg$ has two distinct coefficients at least $m(u)$. If either is larger than $m(u)$ then the claim follows immediately. Otherwise, write their contributions as $m(u)q^i$ and $m(u)q^j$, with $i \in \nn_0$ and $j \in \nn_0$ satisfying $i < j$. Repeatedly using the identity $a q^k = b q^{k + 1}$ for $k \in \nn_0$ gives
	\[
	m(u)q^i + m(u)q^j = (m(u) - a)q^i + (b - a)\sum_{k = i + 1}^{j - 1}q^k + (m(u) + b)q^j.
	\]
	The sum is empty when $j = i + 1$. Every coefficient on the right is nonnegative. Substituting this expression into $fg$ therefore gives a representation of $v$ with a coefficient at least $m(u) + b$. Hence $m(v) > m(u)$. \hfill $\blacksquare$
	
	Suppose by way of contradiction that $r$ does not satisfy the ACCP. Passing to a subsequence of a nonstabilizing chain, we obtain a strictly ascending chain of principal ideals
	\[
	r_0S_q^* \subsetneq r_1S_q^* \subsetneq r_2S_q^* \subsetneq \cdots
	\]
	where $r_i \in S_q^*$ for every $i \in \nn_0$ and $r_0 = r$. Each $r_i$ divides $r$, so Claim~1 shows that $m(r_i)$ exists. Moreover, $m(r_i) \ge a$ for every $i$: otherwise, Lemma~\ref{lem:unique-polynomial-representation} would give a unique representation of $r_i$, and Proposition~\ref{prop:unique-divisors} would give the ACCP for $r_i$, contradicting the infinite strict tail starting there.
	
	Define
	\[
	h_i := \max\{k \in \nn_0 : q^k \mid_{S_q^*} r_i\}
	\]
	for every $i \in \nn_0$, with $h_i \in \nn_0$. This maximum exists, since $\omega_q(q) > 0$ and $k\omega_q(q) \le \omega_q(r_i)$ whenever $q^k$ divides $r_i$. The sequence $(h_i)_{i \ge 0}$ is nonincreasing. If $r_i/r_{i + 1} = q^e$ for some $e \in \nn$ then $h_i \ge h_{i + 1} + e$. Consequently, only finitely many quotients $r_i/r_{i + 1}$ are powers of $q$. Choose $N \in \nn_0$ beyond all those indices, taking $N = 0$ if there are none. Claim~2 then gives
	\[
	m(r_N) > m(r_{N + 1}) > m(r_{N + 2}) > \cdots,
	\]
	an impossible strictly decreasing sequence of nonnegative integers. Therefore $r$ satisfies the ACCP.
\end{proof}

The following example distinguishes the number of representations from the size of their coefficients and explains why the restriction $q < 1/2$ in Lemma~\ref{lem:accp-unique-representation} is necessary.

\begin{example}
	Set $q := 5/6 \in \qq$ and $r := 5 \in S_q^*$. The polynomial representations of $r$ are exactly the constant polynomial $5$ and the polynomials $f_n \in \nn_0[x]$ defined by
	\[
	f_n(x) := x + \cdots + x^{n - 1} + 6x^n
	\]
	for $n \in \nn$, where the sum is empty when $n = 1$. Repeatedly using the identity $5q^i = 6q^{i + 1}$ for $i \in \nn_0$ verifies that $f_n(q) = 5$. To prove completeness, let $f \in \nn_0[x]$ satisfy $f(q) = 5$, and write $f - 5 = (6x - 5)h$ with $h \in \zz[x]$. Set $h_i := [x^i]h \in \zz$ for $i \in \nn_0$. The sequence $(h_i)_{i \ge 0}$ is eventually zero. If $h_i < 0$ for some $i$ then nonnegativity of $[x^{i + 1}]f = 6h_i - 5h_{i + 1}$ forces $h_{i + 1} < 0$, contradicting the eventual vanishing of the coefficients. Thus, each $h_i \ge 0$. The constant coefficient of $f$ gives $h_0 \le 1$, while $5h_{i + 1} \le 6h_i$ implies that the nonzero coefficients of $h$, when present, form an initial block of $1$s. This gives precisely the listed representations. Their degrees are unbounded, but their coefficients are at most $6$. Since $B_q = 24$, Theorem~\ref{thm:accp-coefficient-bound} shows that $5$ satisfies the ACCP. \hfill $\blacksquare$
\end{example}

The next example shows that the coefficient bound can be attained.

\begin{example}
	Let $q := 5/6 \in \qq$, so $B_q = 24$. The element $24$ satisfies the ACCP, whereas $25$ does not. To verify the first assertion, let $f \in \nn_0[x]$ represent $24$, and write $f - 24 = (6x - 5)h$ with $h \in \zz[x]$. As in the preceding example, the coefficients $h_i := [x^i]h \in \zz$, indexed by $i \in \nn_0$ and extended by zero, are nonnegative. The constant coefficient gives $h_0 \le 4$, and the inequalities $5h_{i + 1} \le 6h_i$ imply inductively that every $h_i \le 4$. Hence the constant coefficient of $f$ is at most $24$, and
	\[
	[x^{i + 1}]f = 6h_i - 5h_{i + 1} \le 24
	\]
	for every $i \in \nn_0$. The coefficient criterion therefore gives the ACCP for $24$. In contrast, the constant representation of $25$ has a coefficient greater than $B_q$, so the same theorem rules out the ACCP for $25$. Thus, equality in the coefficient bound is allowed and is attained. \hfill $\blacksquare$
\end{example}

For $q < 1/2$, Lemma~\ref{lem:accp-unique-representation} identifies ACCP elements by uniqueness of representation. Throughout the non-ACCP range, Theorem~\ref{thm:accp-coefficient-bound} identifies them by a uniform coefficient bound. The preceding examples show that bounded coefficients can coexist with infinitely many representations. Thus, the coefficient criterion applies even when uniqueness of representation fails.

\bigskip
\section*{Acknowledgments}
\label{sec:acknowledgments}

This collaboration took place as part of the PRIMES-USA program at MIT. The authors are grateful to the program for making this collaboration possible.


\bigskip
\begin{thebibliography}{99}
	
	\bibitem{ABLST23} K. Ajran, J. Bringas, B. Li, E. Singer, and M. Tirador, \emph{Factorization in additive monoids of evaluation polynomial semirings}, Comm. Algebra \textbf{51} (2023) 4347--4362.
	
	
	
	
	\bibitem{AG22} D.~F. Anderson and F. Gotti, \emph{Bounded and finite factorization domains}, in: Rings, Monoids and Module Theory (A. Badawi and J. Coykendall, eds.), Springer Proc. Math. Stat. \textbf{382}, Springer, Singapore, 2022, pp. 7--57.
	
	\bibitem{AZ07} D.~D. Anderson and M. Zafrullah, \emph{The Schreier property and Gauss' lemma}, Boll. Unione Mat. Ital. (8) \textbf{10-B} (2007) 43--62.
	
	\bibitem{gA22} G. Angerm\"uller, \emph{Strong atoms in monadically Krull monoids}, Semigroup Forum \textbf{104} (2022) 10--17.
	
	\bibitem{BCG21} N.~R. Baeth, S.~T. Chapman, and F. Gotti, \emph{Bi-atomic classes of positive semirings}, Semigroup Forum \textbf{103} (2021) 1--23.
	
	\bibitem{BG20} N.~R. Baeth and F. Gotti, \emph{Factorizations in upper triangular matrices over information semialgebras}, J. Algebra \textbf{562} (2020) 466--496.
	
	\bibitem{BGHL26} G. Blitz, F. Gotti, D. Han, and H. Liang, \emph{Maximal common divisors in monoid algebras}. Preprint, 2026. \href{https://arxiv.org/abs/2609.18947}{arXiv:2609.18947}.
	
	\bibitem{CF19} F. Campanini and A. Facchini, \emph{Factorizations of polynomials with integral non-negative coefficients}, Semigroup Forum \textbf{99} (2019) 317--332.
	
	\bibitem{CCMS09} P. Cesarz, S.~T. Chapman, S. McAdam, and G.~J. Schaeffer: \emph{Elastic properties of some semirings defined by positive systems}, in Commutative Algebra and Its Applications (Eds. M. Fontana, S.~E. Kabbaj, B. Olberding, and I. Swanson), Proceedings of the Fifth International Fez Conference on Commutative Algebra and its Applications, Walter de Gruyter, Berlin, 2009, pp. 89--101.
	
	\bibitem{CGG20} S.~T. Chapman, F. Gotti, and M. Gotti, \emph{Factorization invariants of Puiseux monoids generated by geometric sequences}, Comm. Algebra \textbf{48} (2020) 380--396.
	
	\bibitem{CGG21} S.~T. Chapman, F. Gotti, and M. Gotti, \emph{When is a Puiseux monoid atomic?}, Amer. Math. Monthly \textbf{128} (2021) 302--321.
	
	\bibitem{CGLY26} T. Chen, F. Gotti, T. Lu, and A. Yao, \emph{On the additive structure of algebraic valuations of polynomial semirings II}. Preprint, 2026. \href{https://arxiv.org/abs/2607.02874}{arXiv:2607.02874}.
	
	\bibitem{CGY26} T. Chen, F. Gotti, and A. Yao, \emph{Additive factorizations in finitely generated semidomains}. Preprint, 2026. \href{https://arxiv.org/abs/2610.39309}{arXiv:2610.39309}.
	
	\bibitem{pC17} P.~L. Clark, \emph{The Euclidean criterion for irreducibles}, Amer. Math. Monthly \textbf{124} (2017) 198--216.
	
	\bibitem{CG22} J. Correa-Morris and F. Gotti, \emph{On the additive structure of algebraic valuations of polynomial semirings}, J. Pure Appl. Algebra \textbf{226} (2022) 107104.
	
	\bibitem{DDGLPVZ26} J. Dani, A. Deng, M. Gotti, B. Li, A. Paladiya, J. Vulakh, and J. Zeng, \emph{On the set of atoms and strong atoms in additive monoids of cyclic semidomains}. Comm. Algebra (to appear).
	
	\bibitem{DGHLS25} J. Dani, F. Gotti, L. Hong, B. Li, and S. Schlessinger, \emph{On finitary power monoids of linearly orderable monoids}. Preprint, 2025. \href{https://arxiv.org/abs/2501.03407}{arXiv:2501.03407}.
	
	\bibitem{DGLP26} J. Dani, F. Gotti, B. Li, and A. Paladiya, \emph{Factorizations in algebraic monogenic semidomains}. Preprint, 2026. \href{https://arxiv.org/abs/2609.32138}{arXiv:2609.32138}.
	
	\bibitem{DGZ26} A. Deng, F. Gotti, and J. Zeng, \emph{Factorizations in rational monogenic semidomains}. Preprint, 2026. \href{https://arxiv.org/abs/2607.10178}{arXiv:2607.10178}.
	
	\bibitem{FFNSW26} V. Fadinger, S. Frisch, S. Nakato, D. Smertnig, and D. Windisch, \emph{Primes and absolutely or non-absolutely irreducible elements in atomic domains}, Comm. Algebra \textbf{54} (2026) 2631--2648.
	
	\bibitem{GH06} A. Geroldinger and F. Halter-Koch, \emph{Non-Unique Factorizations. Algebraic, Combinatorial and Analytic Theory}, Chapman \& Hall/CRC (2006).
	
	\bibitem{GGP25} V. Gonzalez, F. Gotti, and I. Panpaliya, \emph{On the ascent of almost and quasi-atomicity to monoid semidomains}. Preprint, 2025. \href{https://arxiv.org/abs/2501.04990}{arXiv:2501.04990}.
	
	\bibitem{fG19} F. Gotti, \emph{Increasing positive monoids of ordered fields are FF-monoids}, J. Algebra \textbf{518} (2019) 40--56.
	
	\bibitem{GG18} F. Gotti and M. Gotti, \emph{Atomicity and boundedness of monotone Puiseux monoids}, Semigroup Forum \textbf{96} (2018) 536--552.
	
	\bibitem{GGHL26} F. Gotti, O. Graia, D. Han, and H. Liang, \emph{The Bi-UF Positive Conjecture for quadratic monogenic semirings and related progress}. Preprint, 2026. \href{https://arxiv.org/abs/2607.22669}{arXiv:2607.22669}.
	
	\bibitem{GGP26} F. Gotti, O. Graia, and N. Panse, \emph{The Bi-LF Positive Conjecture Holds}. Preprint, 2026. \href{https://arxiv.org/abs/2610.40688}{arXiv:2610.40688}.
	
	\bibitem{GK24} F. Gotti and U. Krause, \emph{Riemann zeta functions for Krull monoids}, J. Number Theory \textbf{262} (2024) 134--160.
	
	\bibitem{GK26} F. Gotti and U. Krause, \emph{On the infinitude of irreducibles and primes}. Amer. Math. Monthly (to appear), 2026.
	
	\bibitem{GL22} F. Gotti and B. Li, \emph{Divisibility and a weak ascending chain condition on principal ideals}. Preprint, 2022. \href{https://arxiv.org/abs/2212.06213}{arXiv:2212.06213}.
	
	\bibitem{GP23a} F. Gotti and H. Polo, \emph{On the arithmetic of polynomial semidomains}, Forum Math. \textbf{35} (2023) 1179--1197.
	
	\bibitem{GP23b} F. Gotti and H. Polo, \emph{On the subatomicity of polynomial semidomains}, in: Algebraic, Number Theoretic, and Topological Aspects of Ring Theory (J.-L. Chabert, M. Fontana, S. Frisch, S. Glaz, and K. Johnson, eds.), Springer, Cham, 2023, pp. 197--212.
	
	\bibitem{GZ23} F. Gotti and M. Zafrullah, \emph{Integral domains and the IDF property}, J. Algebra \textbf{614} (2023) 564--591.
	
	\bibitem{nL19} N. Lebowitz-Lockard, \emph{On domains with properties weaker than atomicity}, Comm. Algebra \textbf{47} (2019) 1862--1868.
	
	\bibitem{LWZ24} E. Liang, A. Wang, and L. Zhong, \emph{On maximal common divisors in Puiseux monoids}. Preprint, 2024. \href{https://arxiv.org/abs/2410.09251}{arXiv:2410.09251}.
	
	\bibitem{LRZ23} A. Lin, H. Rabinovitz, and Q. Zhang, \emph{The Furstenberg property in Puiseux monoids}. Preprint, 2023. \href{https://arxiv.org/abs/2309.12372}{arXiv:2309.12372}.
	
	\bibitem{LRT24} C. Liu, P. Rodr\'iguez, and M. Tirador, \emph{Subatomicity in rank-2 lattice monoids}, J. Commut. Algebra \textbf{16} (2024).
	
	\bibitem{OP08} D.~C. Ong and V. Ponomarenko, \emph{The Frobenius number of geometric sequences}, Integers \textbf{8} (2008) \#A33.
	
	\bibitem{sZ22} S. Zhu, \emph{Factorizations in evaluation monoids of Laurent semirings}, Comm. Algebra \textbf{50} (2022) 2719--2730.
	
\end{thebibliography}
\end{document}